\documentclass[11pt]{amsart}
\pdfoutput=1

\usepackage{latexsym, graphicx, epsfig, amsmath, amsfonts,amssymb,subfigure}
\usepackage[bookmarks=false]{hyperref}

\usepackage[left=1.25in,right=1.25in,top=1.5in,bottom=1in]{geometry}

\def\e{{\epsilon}}

\def\be{\boldsymbol{\e}}

\def\mb{\mathbf}
\def\R{\mathbb{R}}
\def\u{\mathbf{u}}

\def\w{\mathbf{w}}

\def\W{\mathbf{W}}
\renewcommand{\d}{\mathbf{d}}
\def\D{\mathbf{D}}
\def\U{\mathbf{U}}

\def\H{\mathbf{H}}
\def\J{\mathcal{J}}
\def\E{\mathbb{E}}

\def\I{\mathbf{I}}

\def\K{\mathbf{K}}

\usepackage{color,xspace}

\usepackage{tikz,pgfplots}
\usepgfplotslibrary{external}
\usetikzlibrary{fit}
\usetikzlibrary{backgrounds}
\pgfrealjobname{main}
\long\def\beginpgfgraphicnamed#1#2\endpgfgraphicnamed{\includegraphics{#1}}
\usepackage{ifthen}

\newcommand{\pinv}{\dagger}

\newcommand{\Ve}[1]{\mathbf{#1}}

\usepackage{epstopdf}
\title[Sequential data assimilation]{Sequential data assimilation with multiple nonlinear models and applications to subsurface flow}

\author{Lun Yang}\thanks{Department of Mathematics and System Sciences, Beihang University, Beijing, 100191, China. Email: lun.yang@buaa.edu.cn}
\author{Akil Narayan}\thanks{Department of Mathematics,  and Scientific Computing and Imaging (SCI) Institute, University of Utah, Salt Lake City, UT 84112, USA. Email: akil@sci.utah.edu. A. Narayan was partially funded by AFOSR FA9550-15-1-0467 and DARPA N660011524053.}
\author{Peng Wang}\thanks{Department of Mathematics and System Sciences, Beihang University, Beijing, 100191, China. Email: wang.peng@buaa.edu.cn. L. Yang and P. Wang were partially funded by National Key Research and Development Program of China (Grant No. 2017YFB0701702) and the Recruitment Program of Global Experts.}

\begin{document}

\maketitle

\begin{abstract}
Complex systems are often described with competing models. Such divergence of interpretation on the system may stem from model fidelity, mathematical simplicity, and more generally, our limited knowledge of the underlying processes. Meanwhile, available but limited observations of system state could further complicates one's prediction choices. Over the years, data assimilation techniques, such as the Kalman filter, have become essential tools for improved system estimation by incorporating both models forecast and measurement; but its potential to mitigate the impacts of aforementioned model-form uncertainty has yet to be developed. Based on an earlier study of Multi-model Kalman filter, we propose a novel framework to assimilate multiple models with observation data for nonlinear systems, using extended Kalman filter, ensemble Kalman filter and particle filter, respectively. Through numerical examples of subsurface flow, we demonstrate that the new assimilation framework provides an effective and improved forecast of system behaviour. 
\end{abstract}

%\begin{keywords}
%  Data assimilation, extended Kalman filter, ensemble Kalman filter, particle filter, uncertainty quantification
%\end{keywords}
%
%% REQUIRED
%\begin{AMS}
%  62M20, 86A05  
%\end{AMS}

%\pagestyle{myheadings}
%\thispagestyle{plain}
%\markboth{L. YANG, P. WANG, A. NARAYAN AND D. XIU}
%{MULTI-MODEL UNCERTAINTY QUANTIFICATION}

\section{Introduction}\label{sec:introduction}

Mathematical models are essential tools to understand and predict the behaviour of physical systems. However, one's lack of knowledge renders such models as imperfect approximations of physical reality, and may lead to various interpretations and mathematical representations. To reduce the discrepancies between model forecast and physical truth, one can utilise data measurements as an alternative source of information. However, in practice these empirical measurements are often noisy and scant in number, scope and resolution, which may also introduce competing constitutive relations and uncertainty in parameters. A judicious combination of imperfect models and sparse, noisy data is necessary in modern simulation paradigms.

%Numerical simulations of mathematical models are essential tools for
%predicting the behavior of physical systems. Myriad numerical
%techniques and approximations are used to simulate physical phenomena
%in fluid dynamics, electromagnetics, chemical systems, astrophysics,
%and more. Since all of these simulations involve approximations,
%uncertainty and error are inevitably present in their predictions. To
%complicate matters, a single physical process may be described by
%\textit{multiple} mathematical models and numerical approximations. In
%addition, one may have access to empirical observations of the
%system---noisy and limited in number, scope, and resolution. A natural
%question to ask is how to combine the models and the
%observational data
%% to form a predicted physical state that is `more correct' than than
%% obtained with any of the models individually?
%to predict the physical state with greater fidelity than can be obtained
%with any of the models individually?

Over the years, data assimilation (DA) has become the primary tool to implement such a combination of models and data. Originally developed for linear systems with noisy measurements, the Kalman filter \cite{kalman_new_1960,KalmanBucy61} minimises a quadratic objective and updates prediction by weighing simulation results and the available data at each stage. Further extensions to address nonlinear systems, such as the extended Kalman filter \cite{Gelb74, Jazwinski70}, the ensemble Kalman filter \cite{Evensen_JGR94, evensen_data_2009} and other variants \cite{Anderson_MWR01,AndersonA_MWR99,BurgersLE_MWR98,TippettABHW_MWR03,WhitakerH_MWR02}, have also been proposed and remain a popular research topic. 

%Various techniques for model averaging and data assimilation, all of
%which attempt , have
%received attention in recent years. In the case of a single dynamical
%model with a stream of noisy observations, the Kalman filter
%\cite{kalman_new_1960,KalmanBucy61} is both simple and remarkably effective. The
%assimilation step of the Kalman filter updates the state by weighing
%the model prediction and the data in order to minimize a 
%quadratic objective. This operation can be interpreted in many different
%ways, for instance, as a minimum variance estimator or as a Bayesian
%update. The original Kalman filter was designed for linear systems,
%but derivate methods for filtering nonlinear systems are plentiful,
%e.g., the extended Kalman filter \cite{Gelb74, Jazwinski70}, the ensemble
%Kalman filter \cite{Evensen_JGR94, evensen_data_2009} and its variants
%\cite{Anderson_MWR01,AndersonA_MWR99,BurgersLE_MWR98,TippettABHW_MWR03,WhitakerH_MWR02},
%and is the subject of extensive ongoing work. 
%%We will attempt to conduct
%%a thorough review here.
%

While DA has been intensely developed, much less work on DA has been conducted to address the problem with multiple simulation models. For static measurements, there are Bayesian model averaging (BMA) \cite{hoeting:1999:bma} and non-Bayesian approaches \cite{hoeting:1999:bma, diks:2010:cpf} to compute predictive distribution and point predictions, respectively.  For dynamic settings where data arrive sequentially, one might consider methods including Dynamic Model Averaging (DMA) \cite{raftery:2010:opu}, the Interacting Multiple Model filter \cite{h._a._p._blom_efficient_1984} and the generalised pseudo-Bayes framework \cite{watanabe_generalized_1993}. A limitation of those methods is that their model averaging weights are scalar-valued; and hence the assimilation process is fixed and cannot account for a changing hierarchy among various models. We also note here that the well-known BMA framework represents a convex average of all simulation results. In other words, its result may be no better than the best available model when all models exhibit consistent bias towards one side of the prediction.

To address such limitations, a new DA framework \cite{Narayan-2012-Sequential} has been recently developed for multiple models. By minimising a variational functional similar to that of Kalman filter, the proposed method introduces a general form of all the models and data, which allows different weights in different Sections of the system state vector. Though the new framework employs a linear filter technique, it is not tied to any specific algorithm for covariance propagation and can be readily combined with most variants of the Kalman filter.

In this paper, we aim to extend the aforementioned discrete-time DA method to nonlinear filtering techniques including the extended Kalman filter, the ensemble Kalman filter and the particle filte. In Section \ref{sec:setup} we formulate the assimilation problem with a review of existing method. Section \ref{sec:setup} introduces our multi-model assimilation algorithm using extended Kalman filter, ensemble Kalman filter and particle filter. To investigate their effectiveness, we conduct numerical simulations in the context of subsurface flow in Section \ref{sec:examples}. Key findings of the new frameworks are concluded in Section \ref{sec:summary}.

%In Section \ref{sec:setup} we formalize notation and discuss
%the task of computing matrix harmonic means and parallel sums, which
%are integral to the method. Section \ref{sec:mahalanobis-mean} extends
%these results to finding a minimal-covariance arithmetic mean of
%independent vector-valued random variables; we call this quantity the
%`Mahalanobis mean,' after the distance functional that it minimizes.
%Section \ref{sec:data-assimilation} discusses the application of the
%techniques from Section \ref{sec:preliminaries} to the problem of
%multi-model assimilation and the Kalman filter. Three proof-of-concept
%example applications are discussed in Section \ref{sec:examples},
%followed by a summary of results in Section \ref{sec:conclusion}.
%

\section{Problem setup} \label{sec:setup}

In this Section we formulate the general problem of data assimilation with multiple models (Section \ref{sec:thesetup}), and present the sequential algorithm of previous work \cite{Narayan-2012-Sequential} on linear systems (Section \ref{sec:oldwork}). 

Our work is conducted on a complete probability space $(\Omega, \mathcal{F}, \mu)$, where $\Omega$ is the collection of events, $\mathcal{F}$ represents a $\sigma$-algebra on sets of $\Omega$, and $\mu$ denotes a probability measure on $\mathcal{F}$. We assume that all random variables have finite second moment.

%The $L^2$ stochastic space is prescribed to all random variables, 
%%
%\begin{equation} \label{eqn:l2space}
%	u \in L^2_\mu \Rightarrow
%	\int_\Omega \| u(\omega)\|^2 \dx{\mu(\omega)} = \omeganorm{{u}(\omega)}^2 \equiv \E \eucnorm{{u}}^2 < \infty.
%\end{equation}
%%
%To ensure the existence of the mean and variance of a vector-valued $u$, $\eucnorm{u}$ denotes the standard Euclidean norm. In the context of limits of random variables, we define equality in the
%$L^2_\mu$ sense: 
%%
%\begin{equation} \label{eqn:l2equality}
%	\lim_{\varepsilon \rightarrow 0} {u}_{\varepsilon} ={u} \Rightarrow \omeganorm{{u}_\varepsilon - {u}} \rightarrow 0. 
%\end{equation}
%%

Vectors $\u$ will be represented by lowercase boldface letters, while uppercase boldface letters will be reserved for matrices $\U$. The superscripts in $\U^{\mathrm T}$ and $\U^\pinv$ are the matrix transpose, and the Moore-Penrose pseudoinverse of $\U$, respectively.  The covariance matrices of random vectors will be denoted by the same letter but with uppercase, e.g. $\U$ is the covariance matrix of a random vector $\u$. %Lastly, the space of all $N\times N$ positive definite matrices is defined by $\pdm$, while $\spdm$ corresponds to the space of all positive semi-definite matrices. 

%%%%%%%%%%%%%%%%%%%%%%%%%%%%%%%%%%%%%%%%%%%%%%%%%%%%%%%%%%%
\subsection{Data assimilation with multiple models}
\label{sec:thesetup}

Let $\u_T\in\R^{N_t}$, $N_t\geq 1$ be the ``true state'' of a physical system that we hope to estimate. The system's dynamics are not completely known but can be approximated by a number of mathematical models $g_m[\cdot], \; m = 1 \dots M $, whose discrete form is
\begin{equation} \label{dynamic-ut}
	%\begin{split}
		%&\frac{\mathrm d \u_T}{\mathrm d t} \approx \mathbf g_m
		%	= \frac{\mathrm d\u_m}{\mathrm d t} , \qquad t\in (0,T], \quad T > 0.  \\
		\u_T(t+1)= g_m\left[ \, \u_T(t) \,\right] + \be_m(t) \\
		%&\u_T(0) = \u_{T0}.
	%\end{split}
\end{equation}
where $ \be_m(t)\in\R^{N_m}$ is a zero-mean random field representing each model's error. 

This "true state" is usually unavailable but may be estimated by the model forecast $\u_m(t) \in \R ^{N_m}$:
\begin{equation} \label{forecast-um}
		 \u_m(t+1)  = g_m\left[  \, \u_m (t) \,  \right],
\end{equation}
We model this forecast state as a random variable, whose error compared to the true state is associated to its variance. Such errors are then characterized by a covariance matrix $\U\in\mathbb{R}^{N_m\times N_m} $. 

In addition to model forecast information, a set of limited measurements ($\d\in\R^{N_d}$, $N_d\geq 1$)  may be at one's disposal but are subject to various sources of uncertainty, including but not limited to measurement instrument error and numerical roundoff error. We model this collective data uncertainty as a zero-mean random field $\be_d\in\R^{N_d}$ with covariance matrix $\D=\E[\be_d\be_d^{\mathrm T}] \in\R^{N_d\times N_d}$.
\begin{eqnarray}	\label{d}
	  \d = \mathcal{H}_d(\u_T) + \be_d, \quad \E[\be_d]=\mb{0},
\end{eqnarray}
The operator $\mathcal{H}_d(\cdot)$ is a matter of assimilation choice and can be linear or nonlinear, e.g., \cite{rabier_overview_2005}. For simplicity, we consider the linear case $\mathcal{H}_d = \H = \I$, e.g. the identity matrix. 

With model forecast and data at hand, our goal is to assimilate this information to compute an ``analyzed state" $\w \in \R^{N_t}$ that should closely resemble the ``true state" at an arbitrary time $t$:
\begin{equation}\label{eq:general-assimilation}
	 \u_T(t) \approx \Ve{w}(t)  = \mathcal L \left[ \Ve{u}_m(t) ,\, \d \right], 
\end{equation}
with an error covariance matrix $\W \in\R^{N_t\times N_t}$. The state $\w$ is chosen as that which minimizes the assimilation objective $\mathcal L \left( \cdot \right)$, which we take to be the following standard functional:
\begin{equation} \label{KF_min}
	\J[\w] = (\w-\u)^{\mathrm T} \U^{-1}(\w-\u) + (\H\w - \d)^{\mathrm T} \D^{-1} (\H\w-\d).
\end{equation}
We note that the functional $\J[\w]$ represents the sum of Mahalanobis distances between the known states {$\u$} and {$\d$}. Our goal is to find an intermediate state {$\w$} that is as close to both states. Another interpretation is to view the computation of {$\w$} as a Bayesian update, in which {$\u$} and its covariance would specify a prior Gaussian distribution on state space; likewise we assume that {$\d$} is obtained as a Gaussian perturbation from linear measurements of the truth, the likelihood of observing the data can be then computed; the analyzed state {$\w$} is the mode of the posterior; equivalently, {$\w$} minimizes the negative log-likelihood of the posterior. 

It is noted that the discrepancy between analyzed state and true state at current time will contribute to model forecast error and its covariance at next time step $t+1$
\begin{eqnarray}
	\label{eqn:forecast-covariance}
        \u_T(t+1)- \u_m (t+1)=  g_m \left[\, \u_T(t) \,\right] - g_m \mathbf{\left[ \, \w(t) \,\right]} + \be_m(t).  
\end{eqnarray}
%

%for some linear/nonlinear operators $\mathcal{H}_m\in \R^{N_m\times N_t}$ and measurement matrix $\mathcal{H}_d \in \R^{N_d\times N_t}$.  

%%%%%%%%%%%%%%%%%%%%%%%%%%%%%%%%%%%%%%%%%%%%%%%%%%%%%%%%%%%
\subsection{Iterative assimilation with multiple models}
\label{sec:oldwork}

A recent study \cite{Narayan-2012-Sequential} proposed an effective DA framework for multi-model forecasts and data assimilations. The algorithm can be implemented iteratively for linear models $\mathbf g_m(\cdot)=A_m(\cdot), \, A_m\in \R^N_m$:
\begin{itemize}
	\item [] {\em At simulation timestep $t \ge 1$},

	\begin{itemize}
	\item {\em Initialization.} Use model forecast $\u_1$ and data $\d$ to perform the standard Kalman update, obtain an ``analyzed state" $\w_1$ and its covariance $\W_1$:
		\begin{eqnarray} \label{algorithm:initial}
			\begin{split} 	
				& \mb{K}_1=\U_1\H^{\mathrm T}(\H\U_1 \H^{\mathrm T}+\D)^{-1}, \\
				& \w_1 = \u_1 + \mb{K}_1(\d-\H\u_1), \\
				& \W_1 = (\mb{I}-\mb{K}_1\H)\U_1 (\mb{I}-\mb{K}_1\H)^{\mathrm T}+\mb{K}_1\mb{DK}_1^{\mathrm T}=(\mb{I}-\mb{K_1H})\U_1.
			\end{split}
		\end{eqnarray}

	\item {\em Iteration.} Use the previous analyzed state $\w_1$ as ``model forecast" and new model forecast $\u_2$ as ``data", perform the standard Kalman update, obtain a new ``analyzed state" $\w_2$ and its covariance $\W_2$. Repeat for all models, $m=2, \dots, M$.  The measurement matrices are set as identity matrix $\H_m=\H = \I$:
		\begin{eqnarray}
				\K_m &=& \W_{m-1} \H_m^{\mathrm T} \left( \H_m \W_{m-1} \H_m^{\mathrm T}+ \U_m \right)^\pinv, \nonumber \\
				\label{algorithm:iteration}
				\w_m &=& \w_{m-1} + \K_m (\u_m - \H_m \w_{m-1}) \\
				&=& \w_{m-1} + \W_{m-1} \H_m^{\mathrm T} \left( \H_m \W_{m-1} \H_m^{\mathrm T}+ \U_m \right)^\pinv
						(\u_m-\H_m\w_{m-1}),\nonumber\\  
				\W_m &=& (\I - \K_m \H_m) \W_{m-1} \nonumber\\
				&=& \W_{m-1} - \W_{m-1} \H_m^{\mathrm T} \left(\H_m \W_{m-1} \H_m^{\mathrm T} + \U_m\right)^\pinv \H_m \W_{m-1}. \nonumber 
		\end{eqnarray}

	\item {\em Forecast.} Use \eqref{forecast-um} $\&$ \eqref{eqn:forecast-covariance} to compute the ``forecast state" for all models $\u(t+1)$ and update its covariance matrix $\U(t+1)$ at the next time step:
	\begin{eqnarray}
		\u_m(t+1)  &=& A_m \w_M (t),  \nonumber \\
				\label{eqn:forecast-covariance-kf}
		\U_m(t+1) &=& A_m \; A_m ^{\mathrm T} \; \W_M (t)  +  \E
		\left[\be_m(t) \be_m^{\mathrm T} (t)\right].  
	\end{eqnarray}

	\end{itemize}
%	\item {\em Repeat the procedure above until the desired time of interest}.
\end{itemize}

We noted here that the framework is essentially a sequential application of a standard Kalman filter update for each new model, except that the inverse operator is replaced with pseudoinverse. The final assimilated result $\w_M$ is independent of the model/data ordering, provided that $\u_1, \dots, \u_M$ and $\d$ are ``consistent" random variables \cite{Narayan-2012-Sequential}. %Capable of dealing with semi-positive covariance matrices, 
The iterative scheme allows one to assimilate all sources of available information, e.g. subset of models, data or both.
% as it was proved that it minimizes a cost function similar to that of Kalman filter but combining all simulations. 

\section{Assimilation of multiple models}
\label{sec:method}

In this Section, we present three DA frameworks to assimilate multiple model forecasts and data for nonlinear systems. To be specific, our methods would be based on extended Kalman filter (Section \ref{sec:e-kalman}), ensemble Kalman filter (Section \ref{sec:ens-kalman}) and particle filter (Section \ref{sec:particle}), respectively. Each subSection starts with a brief summary of the scheme and follows with the multi-model algorithm. For detailed derivation and discussion of those data assimilation techniques over a single model, we refer the reader to \cite{evensen_data_2009}.

%%%%%%%%%%%%%%%%%%%%%%%%%%%%%%%%%%%%%%%%%%%%%%%%%%%%%%%%%%%
%\subSection{Multi-models assimilations with deterministic parameters} 
%\label{sec:method-deterministic}
\subsection{Extended Kalman filter}
\label{sec:e-kalman}

The Extended Kalman filter (EKF) is a data assimilation technique for nonlinear systems. The main idea is to linearize a nonlinear model to propagate states and covariances. In contrast to Kalman filter for which linear operator is \textit{a priori}, the EKF requires only the operator's differentiability, e.g. $ g_m' \equiv \mathrm d g_m / \mathrm d x_m$, and uses a first-order Taylor expansion to represent the ``true state forecast": 
\begin{eqnarray} \label{eqn:Taylor-ekf}
	g_m \left[ \u^T(t) \right] \approx 
	g_m\left[ \w(t) \right] +  g_m' \left[  \w(t) \right] \left[ \u^T(t) - \w(t) \right].
\end{eqnarray}

The EKF scheme for multi-model assimilation can be implemented in an identical procedure as the Kalman filter algorithm (Section \ref{sec:oldwork}), except that at the {\em forecast} stage, the forecast covariance matrix $\U_m(t+1)$ \eqref{eqn:forecast-covariance-kf} is computed by substituting \eqref{eqn:Taylor-ekf} into \eqref{eqn:forecast-covariance}:
\begin{eqnarray}
	\label{eqn:forecast-covariance-ekf}
	\U_m(t+1) \approx  g_m' \left[ \w(t) \right]\; g_m'^{\, \mathrm T} \left[ \w(t) \right] \; \W (t)  +  \E \left[\be_m(t) \,  \be_m^{\mathrm T} (t)\right].  
\end{eqnarray}
%

%%%%%%%%%%%%%%%%%%%%%%%%%%%%%%%%%%%%%%%%%%%%%%%%%%%%%%%%%%%%
\subsection{Ensemble Kalman filter}
\label{sec:ens-kalman}

The Ensemble Kalman Filter (EnKF) is essentially a Monte Carlo implementation of the Bayesian update problem for nonlinear systems. By assuming Gaussian distributions of the analyzed state and the model error, EnKF generates a number of realizations for each ``model forecast" and replaces the forecast covariance matrix $\U_m$ with the sample covariance. The EnKF for multi-model assimilation process scheme is:
\begin{itemize}
	\item {\em At simulation timestep $n=0$},
	
	\begin{itemize}
		%\item {\em Assimilation.} Employ the {\em Initialization} and {\em Iteration} steps in Kalman filter scheme (Section \ref{sec:oldwork}) to assimilate all models forecasts and data, in order to obtain the analyzed state $\w|_{t=1}$ and its covariance $\W|_{t=1}$. 
			
		\item {\em Sampling.} Generate $N_a$ realisations of the initial state using a Gaussian distribution $\mathcal N$ with given mean $u^T(0)$ and covariance $\W$. %\textcolor{red}{Their numerical values will converge as simulation proceeds. (\em REFERENCE) } 
		\begin{eqnarray} \label{algorithm:ensemble-EnKF-analyzed}
			\w^i(0) \sim \mathcal N\left( u^T(0), W \right), \qquad i=1, \dots, N_a.
		\end{eqnarray}
	\end{itemize}	
	
	\item {\em At simulation timestep $n > 0$},
	
	\begin{itemize}
          \item {\em Forecast.} Using $\w^i(n-1)$, compute the model forecast along with a realization of the model error,
            	\begin{eqnarray}  \label{eqn:forecast-enkf}
            		 \u_m^i(n)  = g_m\left[ \w^i(n-1) \right] + \be_m^i(n).
		\end{eqnarray}
		This yields an ensemble of model forecasts of size $M \times N_a$, whose mean and covariance matrix are   
            	\begin{eqnarray} \label{eqn:forecast-covariance-enkf}
			\begin{split}
            		& \E \left[ \u_m(n) \right] = \frac{1}{N_a} \sum_{i=1}^{N_a} \u_m^i(n), \\  
    			& \U_m(n) = \E\left[  \left( \, \u_m^i(n) -  \E\left[ \u_m(n) \right] \, \right)  \left( \, \u_m^i(n) -  \E	\left[ \u_m(n) \right] \, \right)^{\mathrm T}\right]. 
			\end{split}
            	\end{eqnarray}
		\item {\em Assimilation.} Assimilate each sample's ``models forecasts" $\left\{ \u_1^i(n), \dots, \u_M^i(n)\right\}$ with data $\d$, using \eqref{algorithm:initial} and \eqref{algorithm:iteration}, obtain $N_a$ realisations of ``analyzed state", $\w^i(n), \,  i=1, \dots, N_a$.  This can be done in parallel for all samples.
	\end{itemize}
	
\end{itemize}
%

%%%%%%%%%%%%%%%%%%%%%%%%%%%%%%%%%%%%%%%%%%%%%%%%%%%%%%%%%%%%

%%%%%%%%%%%%%%%%%%%%%%%%%%%%%%%%%%%%%%%%%%%%%%%%%%%%%%%%%%%
%\subSection{Multi-models assimilations with random parameters} 
%\label{sec:method-random}
\subsection{Particle filter} \label{sec:particle}

Particle filter (PF), also termed ``sequential Monte Carlo", is an assimilation method widely used in signal and image processing, molecular chemistry, bioinformatics, economics and mathematical finance, and many other fields. Similar to ensemble Kalman filter, it employs a set of particles (samples) to represent the distribution of ``analyzed state". By relaxing any assumptions on the model $\mathbf g_m(\cdot)$ and the state distributions, the PF resamples the particles at each time step when data is available.  To assimilate a number of models and data, we iterate the PF methodology with the following scheme:
 
\begin{itemize}
        \item {\em At simulation timestep $n=0$},
	\begin{itemize}
		\item {\em Sampling.} Generate $N_a$ realisations of the analyzed state using prescribed distribution and initial values. 
	\end{itemize}	
	\item {\em At simulation timestep $n>0$},
	\begin{itemize}
		\item {\em Forecast.} Substitute all ``analyzed state" samples to \eqref{eqn:forecast-enkf} and compute the ``model forecasts", $\left\{ \u_m^i(n) | m=1, \dots, M; i=1, \dots, N_a \right\}$.
		
		\item{\em Assimilation.} Compute the intermediate weight of all particles and models, $i=1, \dots, N_a$. 
		 %For model one $m=1$, weight against available data $\d$:
		%
		\begin{subequations}  \label{eqn:weight-PF}
		\begin{eqnarray}
			\omega_1^i  &\propto&  P(\d \, | \,  \u_1^i), \\
		%\end{equation}
		%
		%For other models, weight against the mean model forecast $\E [u_m] = \sum_{i=1}^{N_a} u_m^i/N_a$:
				%
		%\begin{equation}
			\omega_m^i &\propto&  P\left(  \sum_{j=1}^{N_a} \frac{u_m^j}{N_a} \, | \,   \u_1^i \right), \quad m=2, \dots, M.
		\end{eqnarray}
		\end{subequations}
		Compute the posterior weight for each particle ($i=1, \dots, N_a$):
		\begin{subequations} \label{eqn:pos-weight-pf}
		\begin{eqnarray}
			\omega_i &=& \omega_1^i \cdot \omega_2^i \cdots \omega_M^i,  \\
			\tilde{\omega}^i &=& \omega_i \Big/ \sum_{j=1}^{N_a} \omega_j.
		\end{eqnarray}
		\end{subequations}

              \item{\em Resampling} For model $m$, generate a set of $N_a$ particles from $\left\{\u_1^i(n)\right\}_{i=1}^{N_a}$ according to the posterior likelihood $\left\{\tilde{\omega}^i , \, i=1, \dots, N_a \right\} $. Do this for each $m = 1, \ldots, M$.
		
	\end{itemize}
\end{itemize}
We note that particle filters have many practical limitations \cite{law_data_2015}. While they allow the simulation and data assimilation for nonlinear systems using non-Gaussian distributions, they can require more computational effort than standard Kalman filters and can yield inaccurate results with nearly noise-free models. The methodology above is a standard Sequential Importance Resampling algorithm, and does not address limitations and weaknesses of particle filters. In practice, one would want to use more sophisticated versions of particle filters that ameliorate these concerns.

\begin{figure}[htbp]
\begin{center}
\includegraphics[width=0.49 \textwidth]{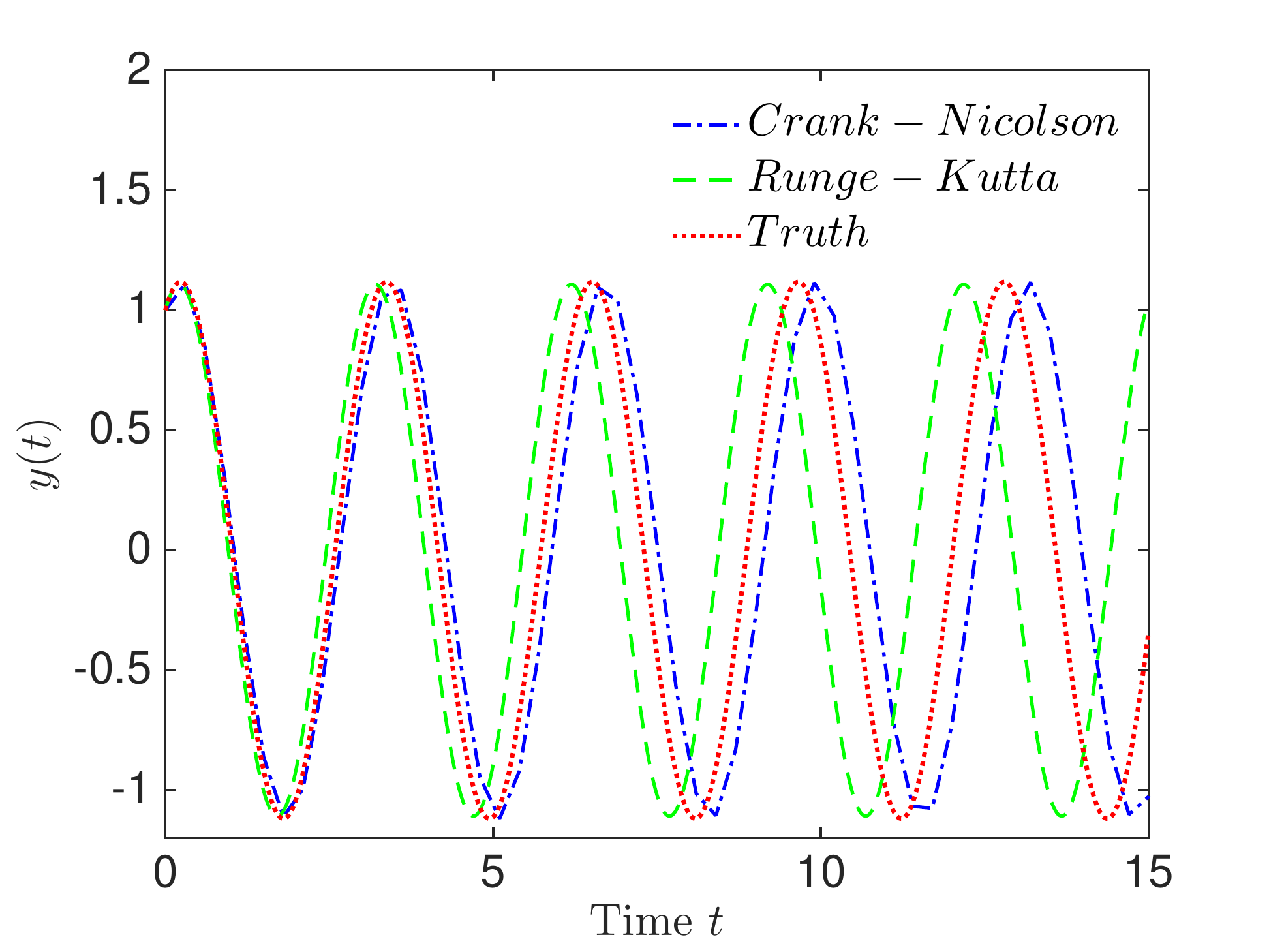} 
\includegraphics[width=0.49 \textwidth]{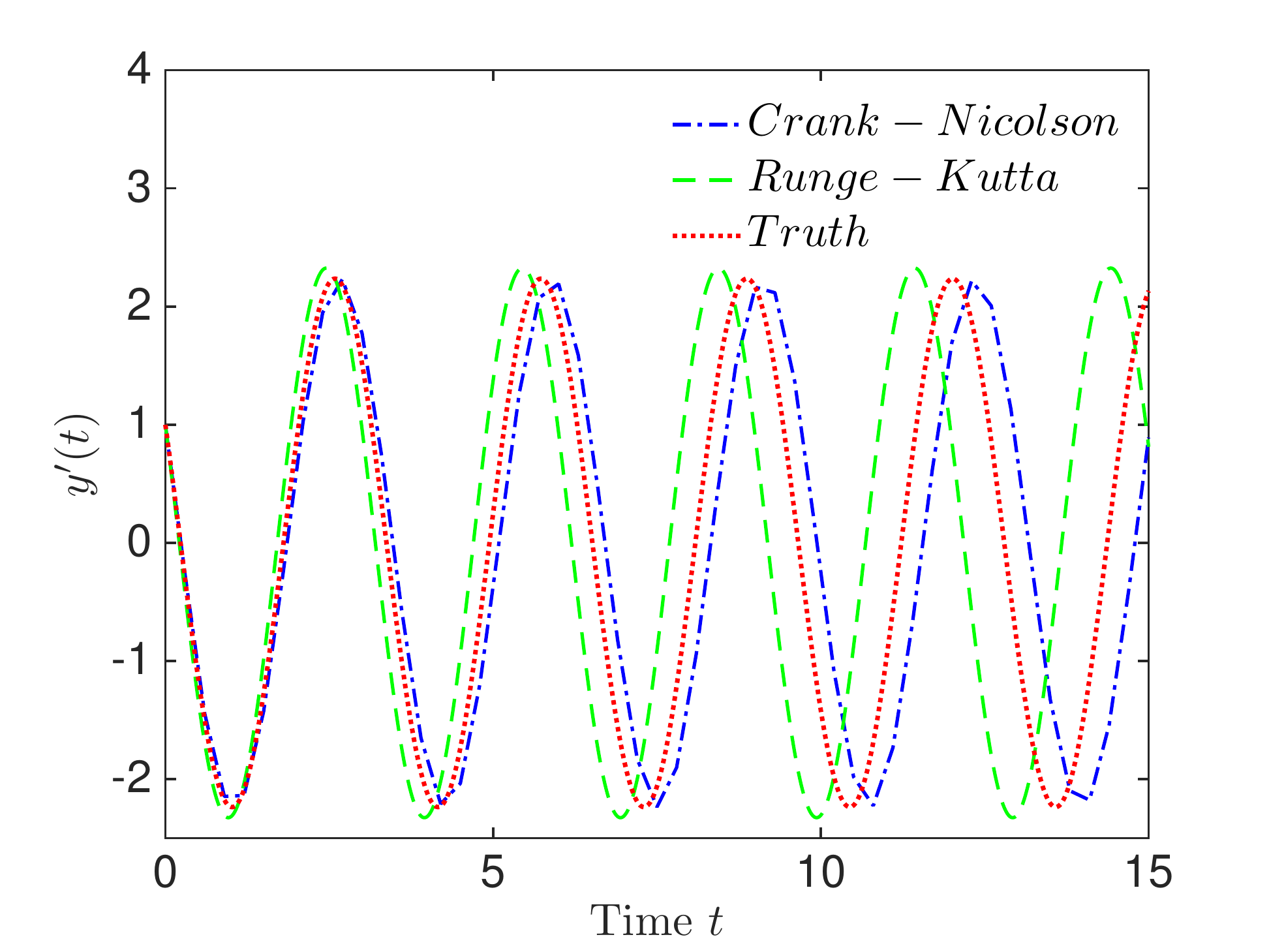} 
\includegraphics[width=0.49 \textwidth]{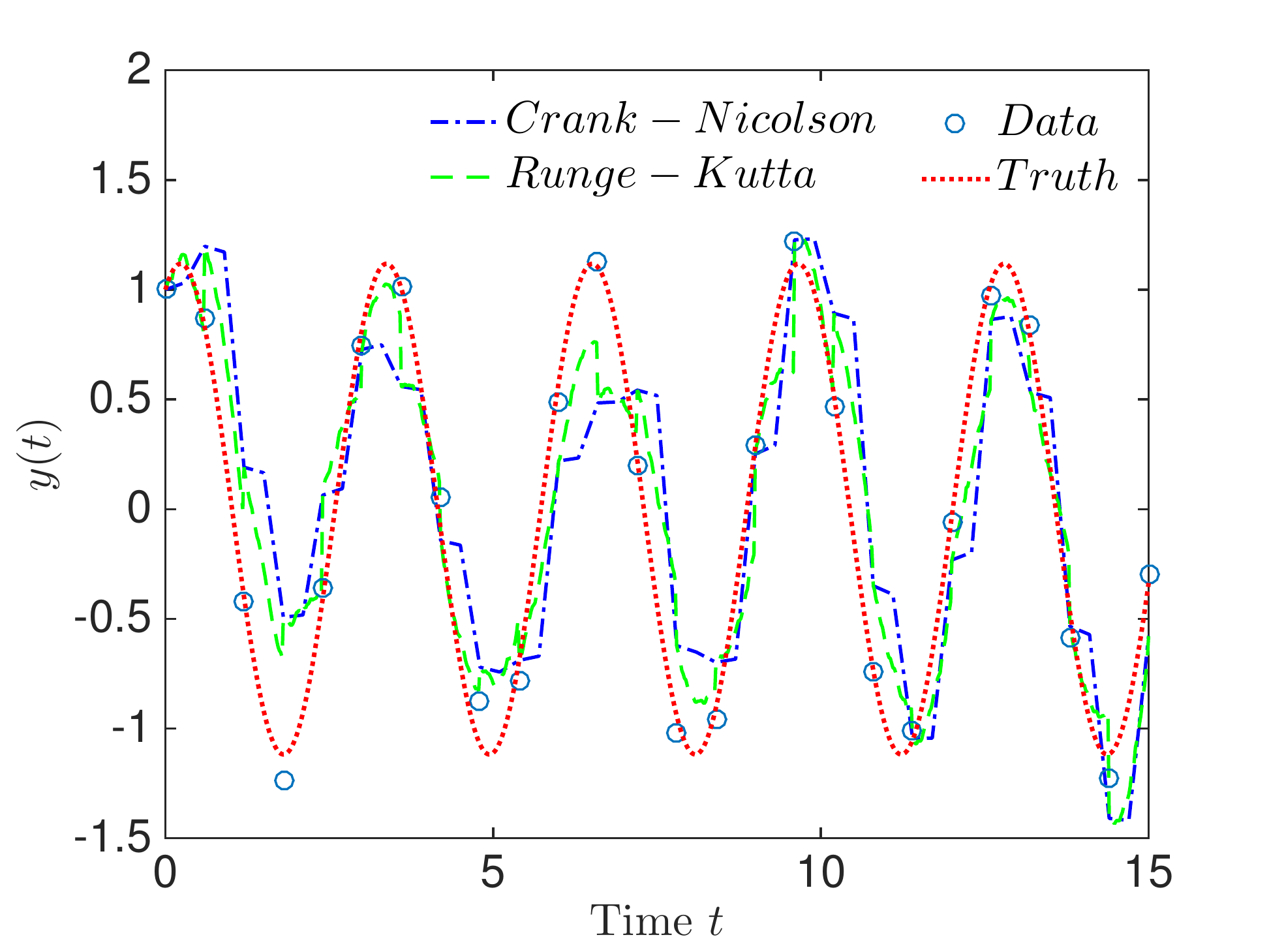} 
\includegraphics[width=0.49 \textwidth]{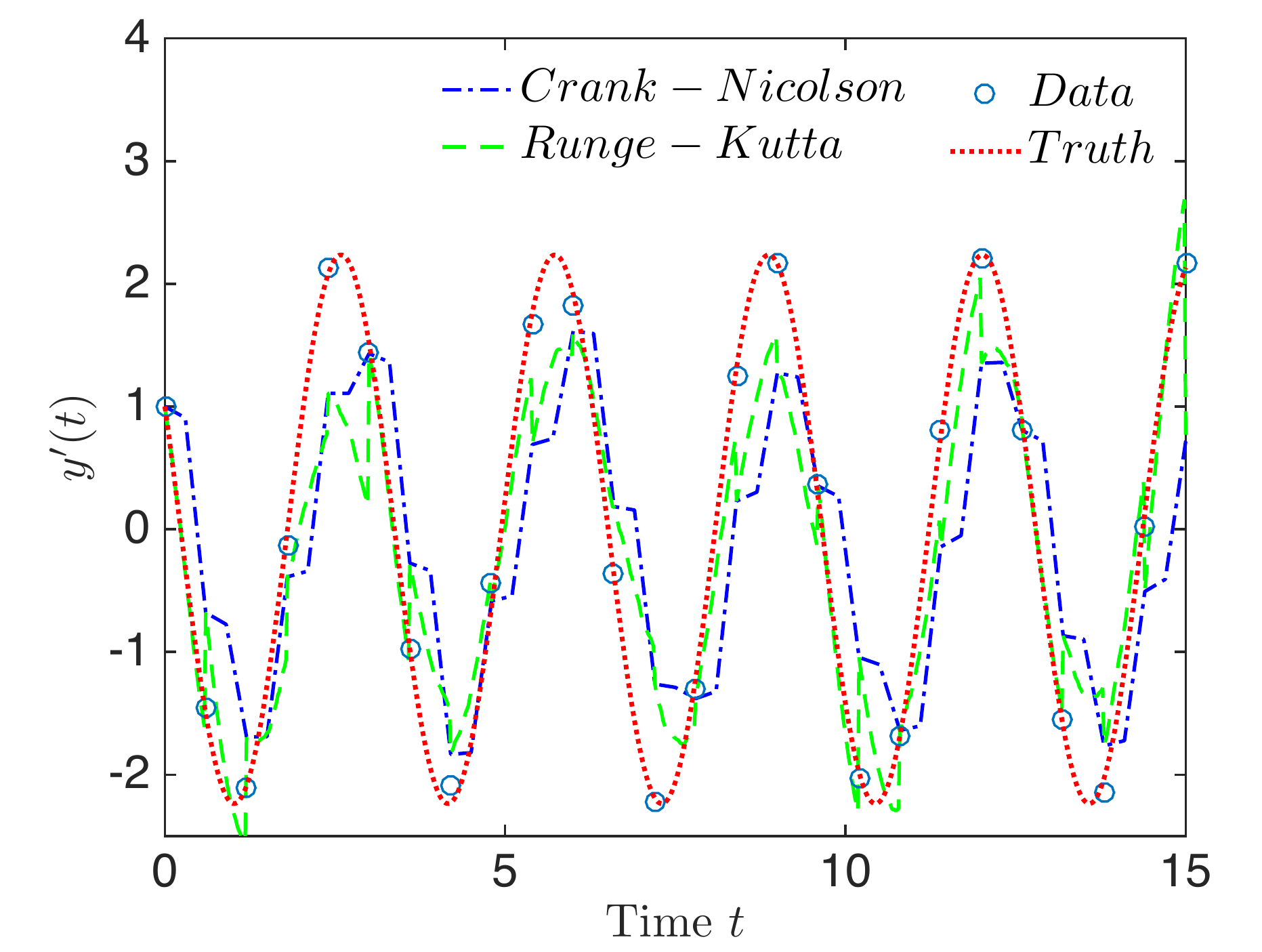} 
\includegraphics[width=0.49 \textwidth]{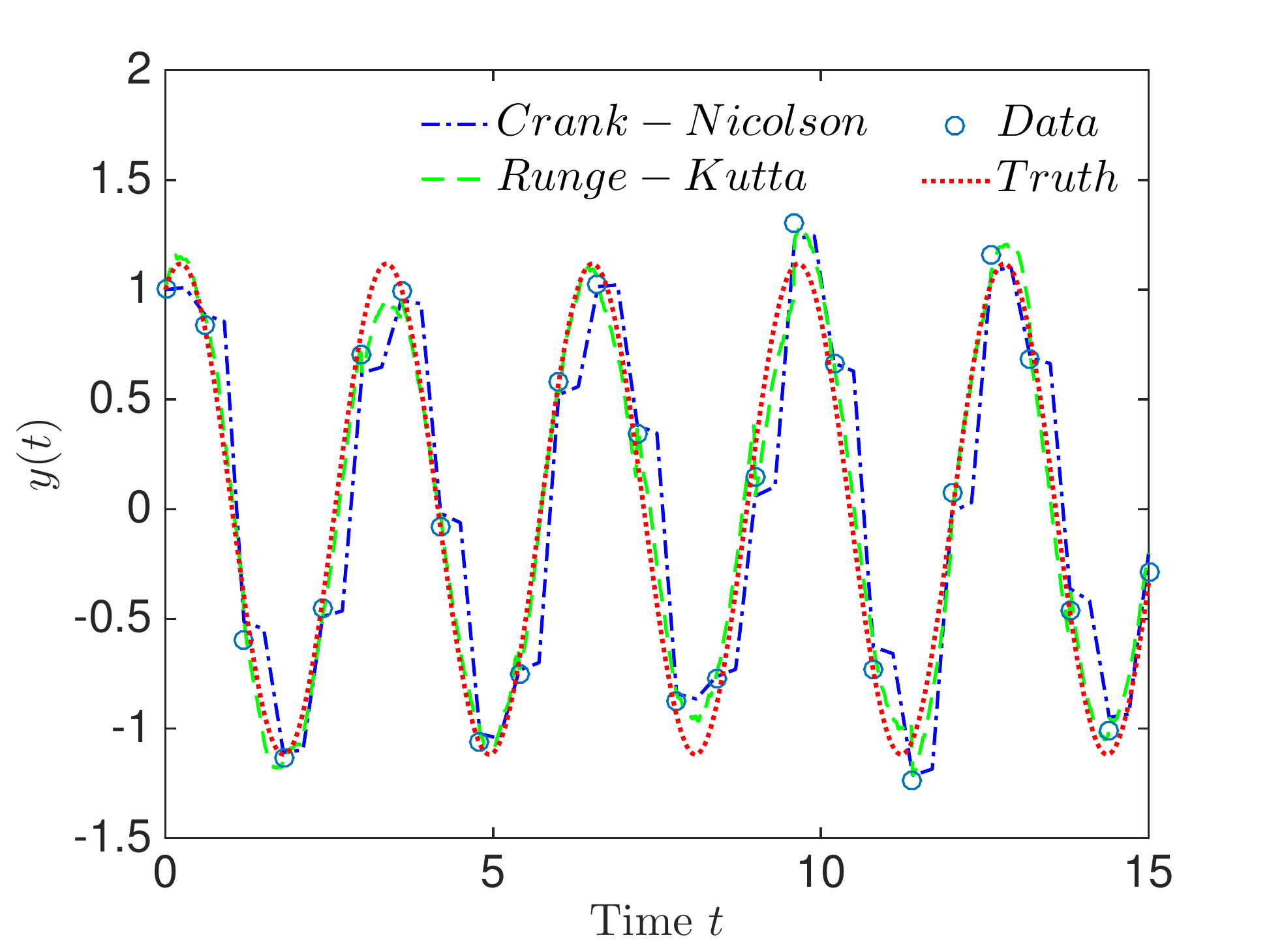} 
\includegraphics[width=0.49 \textwidth]{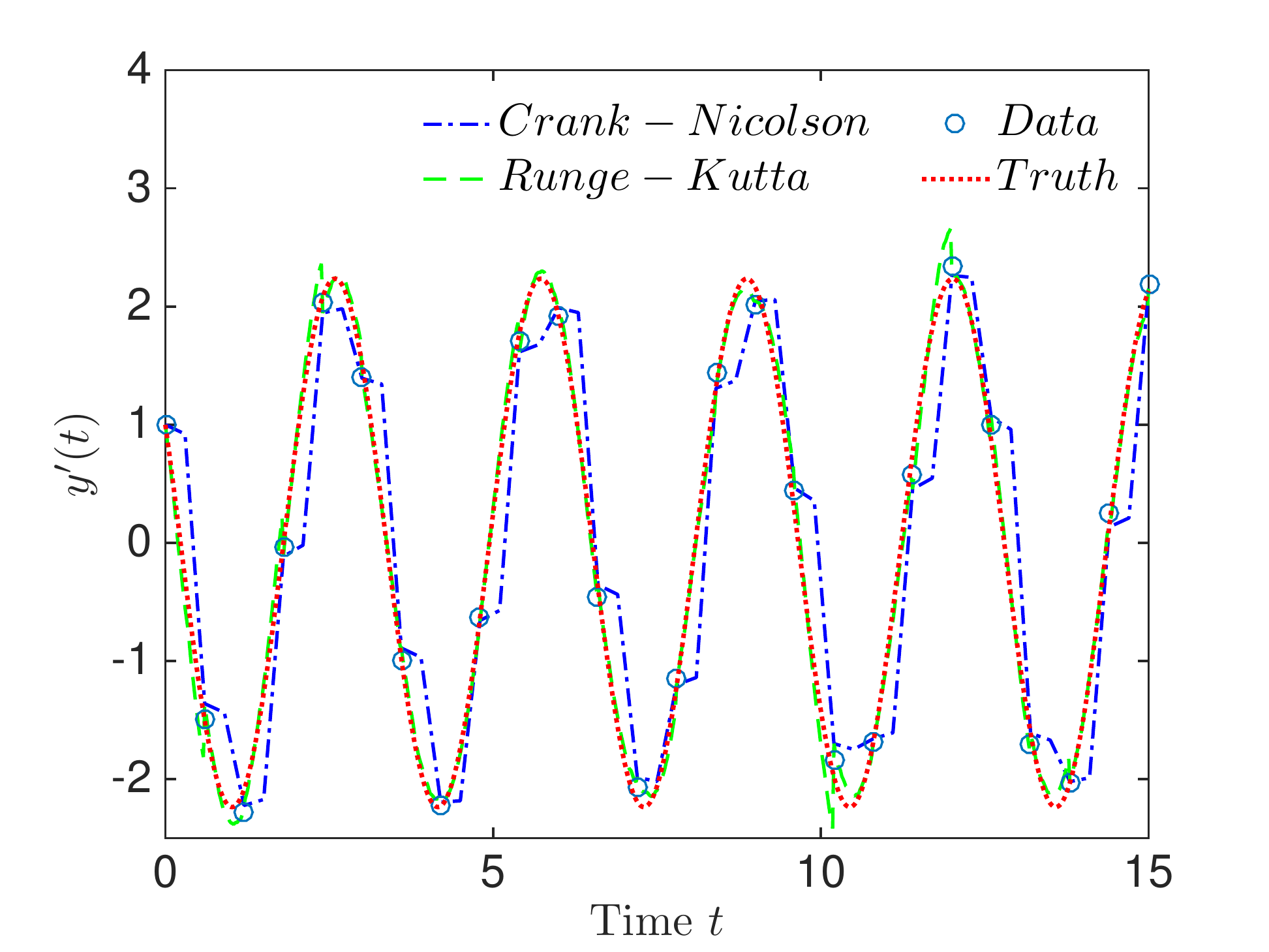} 
  \caption{The solution $(y, y')$ of the system \eqref{eqn:ode} computed using Crank-Nicolson ($CN$) (dashed-dot blue line) and Runge-Kutta of fourth-order ($RK4$) schemes (dashed green line) versus the true solution \eqref{sol:ODE} (dotted red line). Top row: $(y,y')$ without any data assimilation. Center row: $(y,y')$ using the particle filter procedure of Section \ref{sec:particle}, using the Crank-Nicolson model as the reference model (i.e., as model 1). Bottom row: $(y,y')$ using the particle filter procedure of Section \ref{sec:particle}, but with Runge-Kutta as the reference model.} 
\label{fig:pf-ODE} 
\end{center}
\end{figure}

{
%\color{red}
To conclude this section we give a brief motivation for this algorithm: the algorithm above clearly treats model 1 uniquely compared to the other models. We will refer to model 1 as the ``reference model". Assume that one can identify the reference model as the most accurate model we have at our disposal; this determination is application-specific and will depend on a practitioner's expertise. In this case we consider the ensemble $\left\{\u_1^i\right\}_{i=1}^{N_a}$ as the most faithful ensemble of the true, unknown state. We wish to update the distribution of this ensemble using the available data along with the ensembles of the other models. This is acomplished by \eqref{eqn:weight-PF}, where the weights $\omega_m^i$ are updated by comparing the reference model ensemble against the available data and against the mean state of the other models. The formulation \eqref{eqn:pos-weight-pf} combines all these weights into a single weight assigned to each particle in the reference model ensemble. We then resample (as in a standard particle filter strategy) the reference model's ensemble. 

This procedure thus incorporates all model information as well as the data into an ensemble update for each model, using the reference model as the ensemble template. This asymmetric treatment of models has advantages and disadvantages: The advantage is that when a user has a clear preference for one model over all the others, then the preferred model's state is the dominant behavior of the system and all models inherit the accuracy of the preferred model. The disadvantage is that when a user does not have such preferences, then it appears somewhat arbitrary to choose a reference model. 

Our results show that when the reference model is chosen poorly, then the assimilated results are adversely affected. (See Section \ref{sec:shm} and Figure \ref{fig:pf-ODE}.) However, when models have comparable accuracy with no clear preference among them, our results also suggest that the ordering of models in the filtering strategy above affects the results in modest ways. (See Figure \ref{fig:three-scheme}.)

One hopes for a possible convergence statement as the number of particles is increased. We are at present unable to deliver this statement, but hope to provide one in future work.
}

%particle filters
%\cite{van_der_merwe_unscented_2001,leeuwen_nonlinear_2010}

%\textcolor{red}
{
\subsection{Simple example: harmonic oscillator}\label{sec:shm}
We consider a simple example that illustrates some salient properties of the multi-model filtering strategy that we propose above. Consider the following ordinary differential equation modeling a harmonic oscillator
\begin{eqnarray}\label{eqn:ode}
	\begin{split}
	& y'' + w^2y = 0,\\
	& y(0) = y_0,\quad y'(0)=y'_0,
	\end{split}
\end{eqnarray}
with closed-form solution
\begin{eqnarray} \label{sol:ODE}	
  \left(\begin{array}{c}
	y\\
	y'
  \end{array}\right)
	=\left( \begin{matrix} y_0 \cos(wt)+ \frac{y'_0}{w}sin(wt)\\-wy_0\sin(wt)+y'_0\cos(wt)\end{matrix}\right).
\end{eqnarray}
}
%{\color{red}{
We declare the true solution as the function \eqref{sol:ODE} with parameters $w = 2$, $y_0 = 1$, and $y_0' = 1$.  We generate two different approximate models for predicting this solution. The first model solves \eqref{eqn:ode} with $w$, $y_0$, and $y_0'$ as described previously, using a Crank-Nicolson time integrator, where the solution is polluted with Gaussian noise of variance $0.1$. This is an implicit second-order algorithm; because the method is implicit large time steps may be taken, but the resulting solution is relatively inaccurate. Our second model solves \eqref{eqn:ode} with an explicit fourth-order Runge-Kutta method, where $y_0$ and $y_0'$ are as described previously, but instead $w = 2.1$. We again pollute the numerical solution with mean-zero Gaussian noise of variance $0.1$ at every time step. Roughly speaking, the first model predicts the correct system inaccurately, and the second model accurately predicts a system with different parameters. We generate noisy data at specified time points by sampling \eqref{sol:ODE} and then polluting measurements with iid mean-zero Gaussian variables with variance $0.01$.
  
For a fixed $T > 0$ it is in general more expensive to approximate $y(T)$ using the Runge-Kutta method than the Crank-Nicolson method. This is because advancing one time-step with Runge-Kutta is more expensive than advancing one time-step with Crank-Nicolson for the equation \eqref{eqn:ode}, and also because we may take a larger time step with Crank-Nicolson. We set the time step for the Runge-Kutta scheme to $\Delta t = 0.02$, and the time-step for the Crank-Nicolson scheme to $\Delta t = 0.3$. We are provided with noisy data every $0.6$ time units. Thus, data is provided every 2 timesteps of Crank-Nicolson, and every 30 timesteps of Runge-Kutta.
 
As shown in Fig.~\ref{fig:pf-ODE} (top), the two models without data assimilation do not accurately predict the true state. If we set the Crank-Nicolson model as the ``reference" model (i.e., as model 1) and use the multi-model particle filter described in Section \ref{sec:particle}, we obtain results shown in the middle row of Figure \ref{fig:pf-ODE}. The bottom row of Figure \ref{fig:pf-ODE} shows particle filter results using the Runge-Kutta model as the reference model. Both particle filter results become more accurate in terms of predicting frequency, but clearly using Runge-Kutta as the reference model produces superior results. We can understand this behavior from the algorithm description in Section \ref{sec:particle}: while all models are used to update particle weights, only the reference model is used to resample, and therefore the accuracy of the reference model plays an important role in dictating the accuracy of the overall procedure. 

This simple example thus illustrates the following guiding principles of our proposed particle filter:
\begin{itemize}
  \item Model order matters. For our particle filter strategy, this is true even for linear models with Gaussian noise because we are using a nonlinear assimilation technique. In the context of linear models with Gaussian noise, we suspect linear assimilation procedures such as the multi-model Kalman filter \cite{Narayan-2012-Sequential} will prove more attractive: It is known that model ordering for that procedure does not affect the assimilated state. 
  \item An understanding of model accuracy helps. The issue of model order can be ameliorated if a relative estimate of each model's accuracy can be obtained. In particular, choice of the reference model should be made with some knowledge of model accuracy in mind. %The ordering phenomenon observed in Figure \ref{fig:pf-ODE} occurs because the Runge-Kutta integration scheme is more accurate
\end{itemize}
%}
We note that the issues pointed out above are common problems in multi-model simulation strategies. We do not provide constructive solutions to these problems in this manuscript, instead focusing on the problem of assimilation of these models if some prior model ordering or accuracy information is available.
%}
%
%

%
%Compared to other model averaging techniques, the present scheme does
%not need to specify additional dynamics on the model space. Instead,
%the role of the covariance is paramount in determining the weight
%assigned to each model prediction and to the data. The
%covariance-based scheme thus allows matrix-valued weights and is a
%natural generalization of the Kalman filter.

\section{Applications to subsurface flow}\label{sec:examples}

In this section we present the multimodel assimilation framework in the context of subsurface flow into a two-dimensional soil domain $\mathcal{D}$ with $\mathbf x = (x_1, x_3)^{\mathrm{T}} \in \mathcal D$. Details of the three infiltration models are provided in Section \ref{subsec:setting}. We consider a homogeneous  soil domain (Section \ref{subsec:exs-homo}) and examine the impacts of the DA schemes (extended Kalman filter, ensemble Kalman filter and the particle filter) on system prediction (Section \ref{subsec:filter-choice}). This section also includes a study on the time-dependence of model error (Section \ref{subsec:model-error-time}). We then compare our results with that from Bayesian model averaging and investigate their prediction accuracy in the absence of data (Section \ref{subsec:BMA} ). 
%Finally, we introduce parametric uncertainty to our system and apply the particle filter scheme to compute the probabilistic density function (PDF) of infiltration rate $i(t)$, with a reference from Monte Carlo simulations \annote{(Section \ref{subsec:hetero})}.

%provide a few simple examples that illustrate the broad range of applicability of the model assimilation approach previously detailed. Our implemented methods use the iterative procedure outlined in Theorem \ref{thm:iterative-procedure}, and we assume that all models represent consistent random variables. For our examples, this is a valid assumption: all involved random variables have strictly positive covariances. The cost of the assimilation procedure does not suffer greatly from use of the pseudoinverse implementation as a safeguard.     

%Our first example concerns a rudimentary differential equation $y' = a y$ and uses Taylor polynomials as the models. The second example is a related stochastic differential equation (SDE) example that uses the same Taylor polynomial models, but showcases the complex interplay that can occur between the models and the data. Finally, we consider a periodic one-dimensional advection problem with non-constant wavespeed to show how the strengths of different models can be combined by using the multimodel assimilation approach.

\subsection{Infiltration models} \label{subsec:setting}

An infiltration process under ponded condition can be described by several models. Here we consider the Richards equation \cite{Richards-1931-capilary}, the Green-Ampt model \cite{warrick-2003-water}, and the Parlange model \cite{Parlange1982-infiltration}. %All three models consist of the same physical quantities and constitutive relations and all predict the infiltration depth $x_{\mathrm f}$, albeit of various fidelity and computational cost. 
The former is a fine-scale (high-fidelity) model which relies on fewer model assumptions but at a much higher computational cost. The latter two models are coarse-scale and replace the two-dimensional flow with a collection of one-dimensional independent flows (in the $x_3$ direction). Such simplifications is often justified with the Dagan-Bresler parameterization \cite{dagan-1983-unsaturated} and were found to yield accurate predictions in hydrology \cite{bresler-1983-unsaturated,indelman-1998-stochastic,meng-2008-development,morbidelli-2007-simplified,zeller-2000-quantifying,sinsbeck-2015-impact}.

%%%%%%%%%%%%%%%%%%%%%%%%%%%%%%%%%%%%%%%%%%%%%%%%%%
\subsubsection{Richards equation} 
\label{subsec:set-richards}

The Richards equation is a nonlinear advection-diffusion equation that characterizes the movement of water content $\theta(\mathbf{x}, t) : \, \mathcal D \times \mathbb R^+ \rightarrow [\theta_{\mathrm i}, \phi]$ in unsaturated soil:

\begin{eqnarray}\label{eq:richards}
	\frac{\partial \theta}{\partial t} = \nabla \cdot (K \nabla \psi) - \frac{\partial K}{\partial x_3},
	\quad \mathbf x \in \mathcal D = \left\{ -L \le x_1 \le L, \, 0 \le x_3 \le \infty\right\},
	\quad t>0.
\end{eqnarray}

Here $\phi$ denotes soil porosity and $\theta_{\mathrm i}$ is the irreducible water content, whose constant values satisfy $0 < \theta_{\mathrm i} < \phi < 1$. $\psi(\mathbf x, t)$ represents the pressure head and $K(\mathbf x, \theta)$ is the soil hydraulic conductivity. The equation is subject to the initial and boundary conditions:
\begin{subequations}\label{eq:rich-ic-bc}
	\begin{eqnarray}
		&\theta(\mathbf x, 0)= \theta_{\mathrm{init}},
		\qquad
		& \theta(x_1,  x_3=\infty, t) = \theta_{\mathrm{init}}, \\
		 &\psi(x_1,  x_3=0,  t) = \psi_0,
		\qquad 
		&\frac{\partial \psi}{\partial x_3} (x_1 = \pm L, x_3, t) =0.
	\end{eqnarray}
\end{subequations}
The initial water content is no smaller than the irreducible water content, e.g. $\theta_{\mathrm{init}} \ge\theta_{\mathrm i}$.

The infiltration rate $i$ can be obtained as
\begin{eqnarray} \label{eqn:i-xf-richards}
  i(t) = ( \phi - \theta_{\mathrm{init}}) \frac{ \mathrm d x_{\mathrm f} } {\mathrm d t}
	= ( \phi - \theta_{\mathrm{init}}) \frac{\mathrm d}{\mathrm d t} \left[ \int_0^{\infty} \frac{ \theta - \theta_{\mathrm {init}} }{\phi - \theta_{\mathrm {init} }} \mathrm d x_3 \right].
\end{eqnarray}

\paragraph {The constitutive relations} The Richards equation is supplemented by two constitutive relations, $K = K_s K_r(\theta)$ and $\theta = \theta(\psi)$, where $K_s$ and $K_r$ are the saturated and relative hydraulic conductivities, respectively. The constitutive relations reflect soil properties and hence are independent from the choice of flow model. However, since they are derived from data analysis, various interpretations and mathematical representations may exist. As shown in previous section, the DA framework could directly extend to include multiple empirical models of $K_r(\theta)$ and $\theta(\psi)$ in the assimilation process. For simplicity, we restrict the following analysis to a specific constitutive relation, the van Genuchten model \cite{van-1980-closed},
\begin{eqnarray} \label{eqn:vG}
 		K_r= \frac{ [1 - \psi_d^ {mn} (1 + \psi_d^n )^{-m} ]^2 }{ (1+ \psi_d^ n )^{m/2} }, 
	\hskip 25pt
		\frac{\theta - \theta_{\mathrm i}}{\phi -\theta_{\mathrm i}} = \frac{1} {(1+ \psi_d^ n )^{m}},
	\hskip 25pt
		\psi_d = \alpha | \psi |,
	%\quad
	%	m = 1 - \frac{1}{n},
\end{eqnarray}
where $\alpha  > 0$ and $m = 1 - 1/n$ are shape parameters characterising the soil property. 

%%%%%%%%%%%%%%%%%%%%%%%%%%%%%%%%%%%%%%%%%%%%%%%%%%
\subsubsection{Green-Ampt model}
\label{subsec:set-greenampt}
 
The Green-Ampt model \cite{warrick-2003-water} simplifies the two-dimensional flow field to a collection of isolated vertical flow tubes ($x_3$ direction). By assuming the existence of a sharp wetting front $x_{\mathrm f}(t)$ moving downwards, it approximates the S-shaped wetting curve as two separate regions: dry $(\theta = \theta_{\mathrm{init}})$ and wet ($\theta = \phi$), respectively; and arrives at an implicit solution for $x_{\mathrm f}(t)$ and $i(t)$
\begin{eqnarray}\label{xfim}
	&& x_{\mathrm f} - (\psi_0-\psi_{\mathrm f}) \ln \left( 1+ \frac{x_{\mathrm f}}{\psi_0-\psi_{\mathrm f}} \right) = \frac{K_s}{\phi - \theta_{\mathrm{init}} } t, \\
        && \qquad  i(t) = -K_s\frac{\psi_{\mathrm f} - x_{\mathrm f}-\psi_0}{x_{\mathrm f}},
\end{eqnarray}
where $\psi_{\mathrm f}$ is the pressure head at the wetting front and is often set to a capillary drive \cite{Neuman1976_Wetting,  Barry1993-class, wang-2011-reduced}:
\begin{eqnarray} \label{pfbow}
	\psi_{\mathrm f} = -\int\limits_{-\infty}^0 K_r(\psi) \, \mathrm d\psi.
\end{eqnarray}
%
%where $\psi_i=-\infty$ is the pressure head corresponding to the water content $\theta_{\mathrm i}$.  
%%%%%%%%%%%%%%%%%%%%%%%%%%%%%%%%%%%%%%%%%%%%%%%%%%
\subsubsection{Parlange model}
\label{subsec:set-parlange} 

The Parlange et al.~\cite{Parlange1982-infiltration} model approximates the wetting curve $x_{\mathrm f}(t)$ as a sigmoidal form. Under constant ponding water head $\psi_0$ at the surface $(x_3=0)$, it yields an implicit expression for the infiltration rate \cite{Parlange1975-On, Haverkamp1990-infiltration, wang-2011-reduced}:
\begin{subequations}\label{eqn:i-parlange}
	\begin{eqnarray}\label{equ:parlange}
		t &=&\frac{K_s(\psi_0+\psi_{\mathrm j} )(\phi - \theta_{\mathrm{init}} )}{(i-K_s)K_s}
		-\frac{S^2-2\psi_{\mathrm j} K_s (\phi - \theta_{\mathrm{init}} ) }{2 K_s i} \nonumber \\
		&& \qquad +\frac{S^2-2K_s (\phi - \theta_{\mathrm{init}} )(\psi_0+2\psi_{\mathrm j} )}{2(K_s)^2} 
		  \ln(1+\frac{K_s}{i-K_s})
	\end{eqnarray}
Here $\psi_{\mathrm j} $ ($\psi_{\mathrm j} < \psi$) represents a small pressure jump at saturation and remains ``constant in time and independent of changing boundary conditions" \cite{Haverkamp1990-infiltration}. $S$ is the soil sorptivity and takes the following form under a van Genuchten relation~\eqref{eqn:vG}
	\begin{eqnarray}
		S^2=\frac{K_s}{\alpha} (\phi - \theta_{\mathrm{init}}) (1-m)A(m).
	\end{eqnarray}
The constant $A(m)$ is given by
	\begin{eqnarray}
		A(m)&=&\frac{\Gamma(1-m)\Gamma(3m/2-1)}{\Gamma(m/2)}-\frac{4}{3m-2}+\frac{\Gamma(m+1)\Gamma(3m/2-1)}{\Gamma(5m/2)} \nonumber \\
        &&+\frac{\Gamma(1-m)\Gamma(5m/2-1)}{\Gamma(3m/2)}-\frac{4}{5m-2}+\frac{\Gamma(m+1)\Gamma(5m/2-1)}{\Gamma(7m/2)}
	\end{eqnarray}
\end{subequations}
where $\Gamma(\cdot)$ is the Euler Gamma function.
%%%%%%%%%%%%%%%%%%%%%%%%%%%%%%%%%%%%%%%%%%%%%%%%%%
%\subsubsection{Statistics of soil parameters}
%\label{subsec:set-stat}

%%%%%%%%%%%%%%%%%%%%%%%%%%%%%%%%%%%%%%%%%%%%%%%%%%
%\subsubsection{Simulations setup}
%\label{subsec:set-simulation}

%%%%%%%%%%%%%%%%%%%%%%%%%%%%%%%%%%%%%%%%%%%%%%%%%%
\subsection{Simulation setup}
\label{subsec:exs-homo}

In the following simulations, we consider the two-dimensional soil domain as a homogeneous block and set solutions of the Richards equation, computed from the United States Geological Survey simulation code VS2DT \cite{healy-1990-simulations,lappala_documentation_1987}, as the ``true value" of the infiltration rate. The noisy data set is taken from equi-spaced entries within the Richards solution vector, perturbed with noise $\e_d \sim \mathcal N (0,0.002^2)$. It would be assimilated with forecasts from the two reduced-complexity models (Green-Ampt and Parlange). 

\begin{figure}[htbp]
\begin{center}
\includegraphics[width=12 cm]{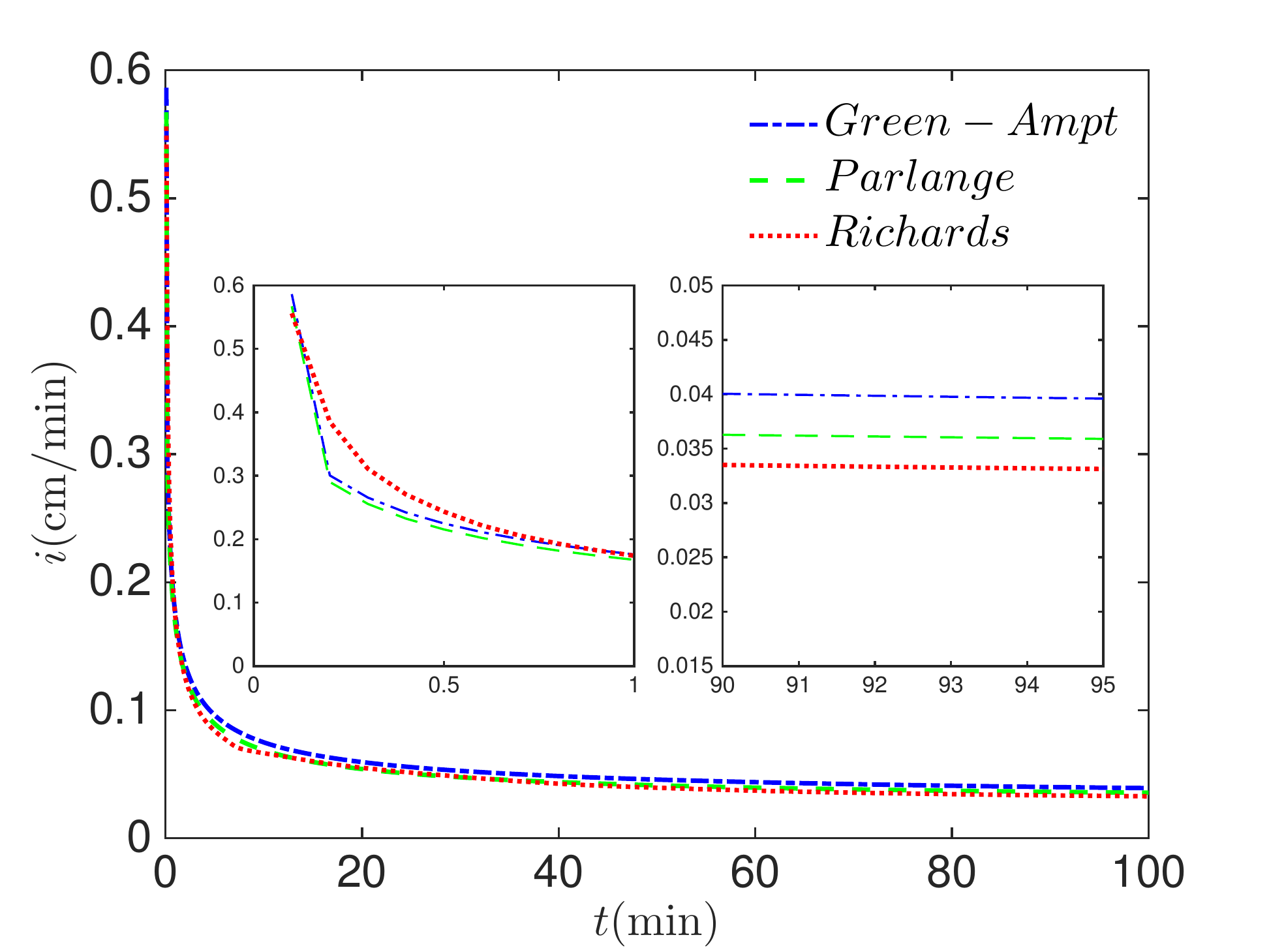} 
\caption{Temporal evolution of infiltration rate $i(t)$ from Green-Ampt model (dash-dot line), Parlange model (dashed line) and Richards equation computed from VS2DT (dotted line) with a close look of the early time and late time.} 
\label{fig:three-models} 
\end{center}
\end{figure}

Figure \ref{fig:three-models} illustrates the differences among Green-Ampt,  Parlange and Richards models by plotting their individual simulation results, without any assimilation. Compared to the Richards solutions, both reduced-complexity models underestimate the infiltration rate at earlier time but overestimate it as time elapses.  %Green-Ampt model exhibits greater error than Parlange model. 
In subsequent study, we take the model errors of Green-Ampt and Parlange model as Gaussian white noises. Each model's error variance is taken as the mean error value between the model and Richards equation over a number of time steps. 

We note that Green-Ampt \eqref{xfim} and Parlange models \eqref{eqn:i-parlange} are both implicit forms for $i(t)$, which can be transformed to dynamical systems by taking a time derivative, with initial infiltration rate $i_0$:
\begin{eqnarray} 
	\label{eqn:initial-greemampt}
	\mbox{Green-Ampt}: \quad 
	\frac{\mathrm d i}{\mathrm d t} &= &-\frac{i(i-K_s)^2}{K_s(\phi - \theta_{\mathrm{init}} )(\psi_0-\psi_{\mathrm f} )}; \\
	\label{eqn:initial-parlange}
	\mbox{Parlange}: \quad
	\frac{\mathrm d i}{\mathrm d t} &=& \frac{2i^2(i-K_s)^2}{S^2(K_s-i)-2K_s (\phi - \theta_{\mathrm{init}}) (\psi_0 i+\psi_{\mathrm j} K_s)}, 
\end{eqnarray}
In the subsequent simulations, a fourth-order Runge-Kutta numerical scheme is applied to solve the two equations above. Unless specified otherwise, Bet-Dagan soil properties \cite{russo-1992-statistical}, as shown in Table ~\ref{pa:table:soil}, are used in subsequent study. 
\begin{table}[hbtp]
\begin{center}
\begin{tabular}{ccccccc}
\hline
 $K_s$ & $\alpha$ & $\phi$  & $\theta_i$ & $\psi_0$  & $\psi_{\mathrm j}$  & van Genuchten $n$ \\
 (cm/min) & (cm$^{-1}$) & & & (cm) & (cm) &  \\
\hline
 $-3.58$ & $-3.01$ & $0.42$ & $0.13$ & $1$ & $2$ & $1.81$\\
%\hline
%variance & $0.89$ & $0.63$ & - & - & - & - & - \\
\hline
\end{tabular}
\caption{Hydraulic properties of the Bet-Dagan soil \cite[Table 3]{russo-1992-statistical}.}
\label{pa:table:soil}
\end{center}
\end{table}

%%%%%%%%%%%%%%%%%%%%%%%%%%%%%%%%%%%%%%%%%%%%%%%%%%
\subsection{Effects of filter scheme}
\label{subsec:filter-choice}
\begin{figure}[htbp]
\begin{center}
\includegraphics[width=0.49\textwidth]{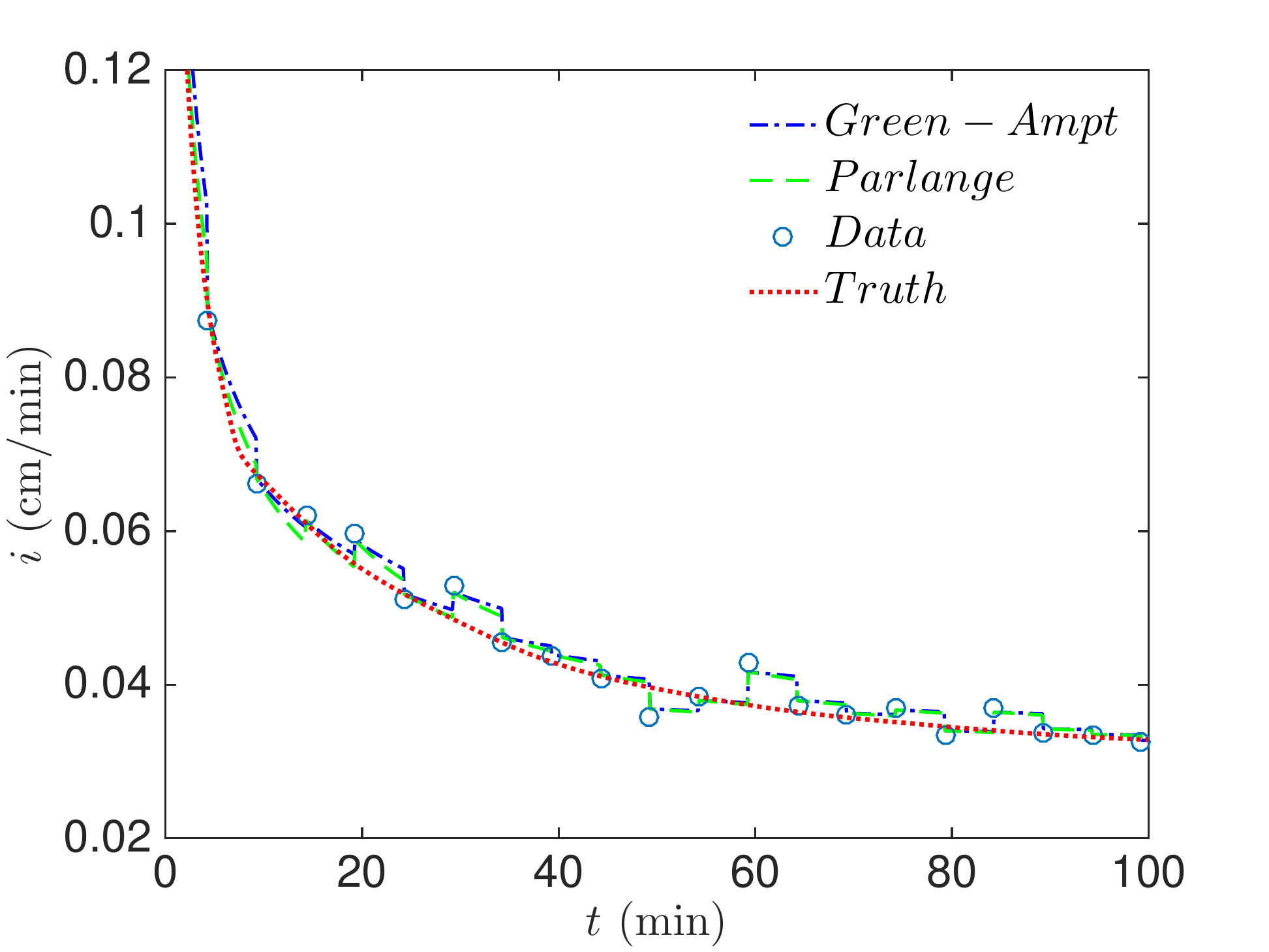} 
\includegraphics[width=0.49\textwidth]{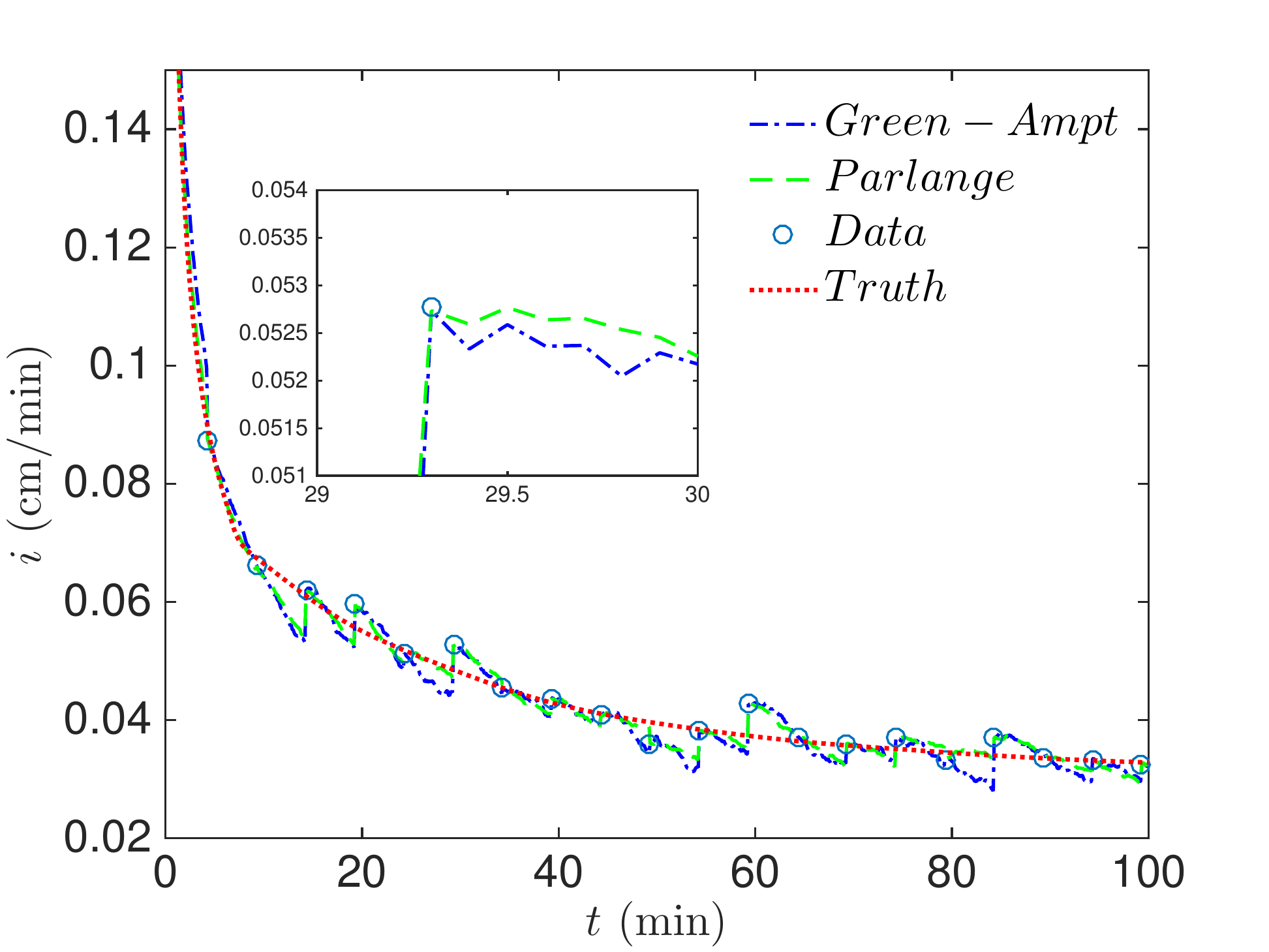}
\includegraphics[width=0.49\textwidth]{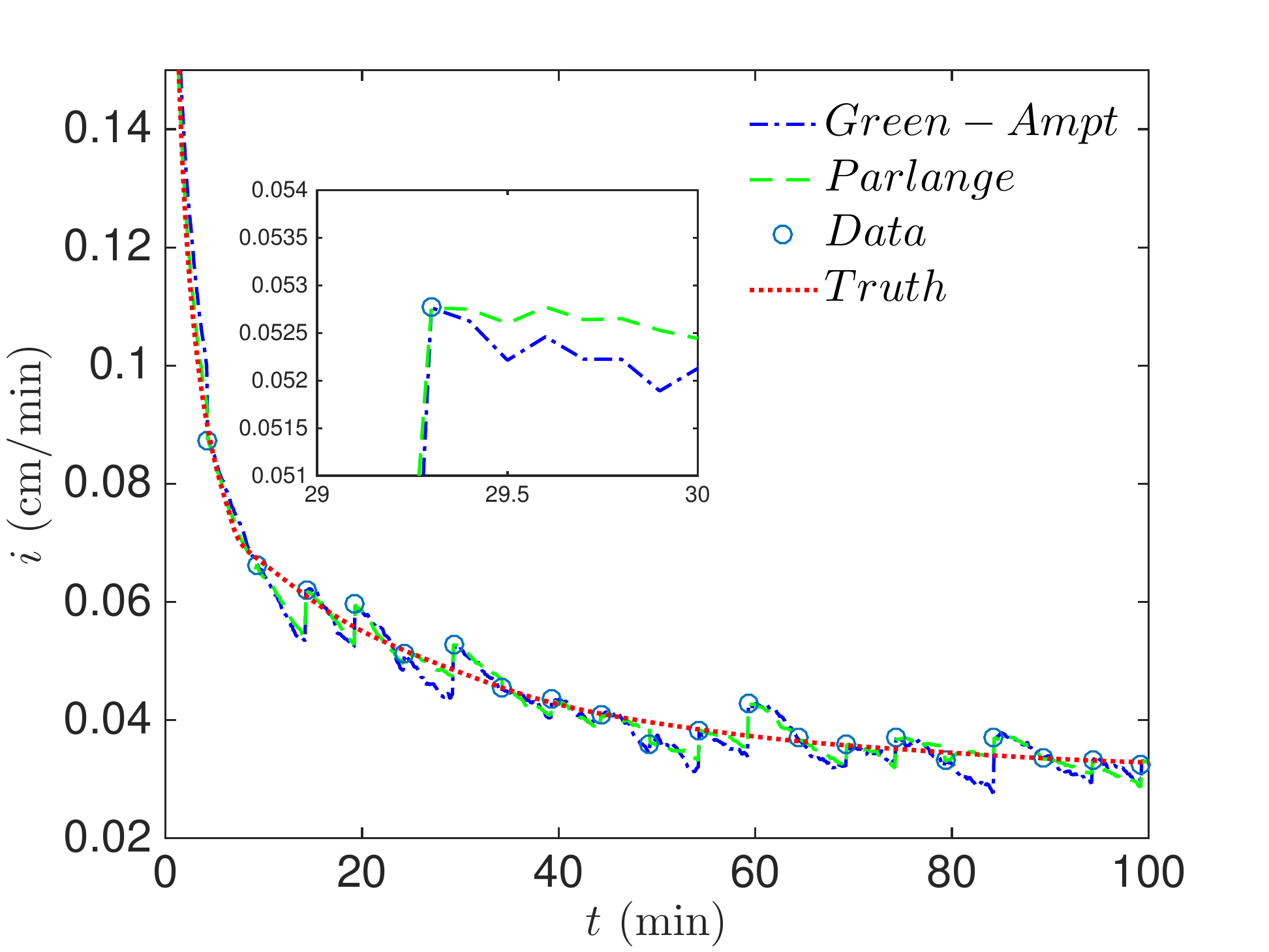}
\includegraphics[width=0.49\textwidth]{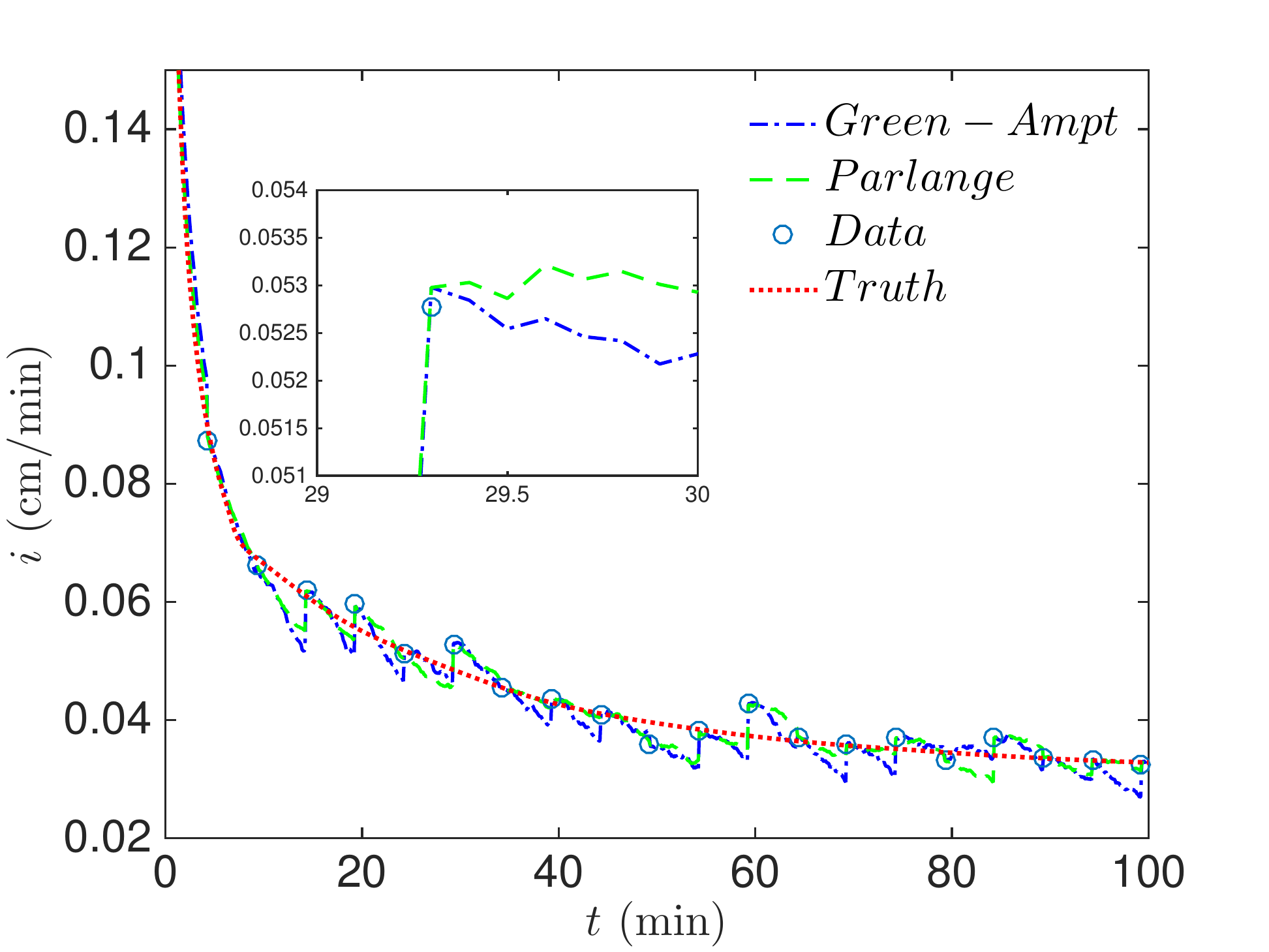} 
  \caption{Temporal evolution of infiltration rate $i(t)$ assimilated from Green-Ampt model (dash-dot line), Parlange model (dashed line) and data (circle), using extended Kalman filter (top left), ensemble Kalman filter with $1000$ samples (top right), a particle filter with $1000$ samples using the Green-Ampt model as the reference model (bottom left), and a particle filter with $1000$ samples using the Parlange model as the reference model (bottom right). The true value computed from VS2DT is denoted as dotted line.} 
\label{fig:three-scheme} 
\end{center}
\end{figure}

Figure \ref{fig:three-scheme} presents the temporal evolution of system states assimilated from the Green-Ampt model, Parlange model, and observation data, using the extended Kalman filter, the ensemble Kalman filter ($1000$ samples) and a particle filter ($1000$ samples), separately. At the assimilation points, all three schemes help reduce discrepancy from individual model forecast or measurement alone.  Between assimilation points, model forecasts vary with the choice of scheme: EKF's model forecasts are relatively parallel to the truth; while the EnKF's and PF's  exhibit sharper temporal change and their profiles are very similar, albeit slight differences still exist at a closer look $t\in[29, 30]$ mins. 

Recalling the discussion of Section \ref{sec:particle} on the multi-model particle filter, we see that the choice of reference model $\u_1^i$  \eqref{eqn:weight-PF} directly affects the computation of the ``analyzed weight" \eqref{eqn:pos-weight-pf} at each assimilation step.  We illustrate the impact of reference model choice $\u^1_i$ on the assimilation results in the bottom-left and bottom-right plots of Fig. \ref{fig:three-scheme}, in which Green-Ampt and Parlange model, respecitvely, is selected as the reference model. Although there are only two models to assimilate, the difference of the reference model choice leads to divergence of predictions: results from the Green-Ampt model occasionally overestimates the infiltration rate, whereas the one from Parlange model fits the true values more closely. One may conclude that Green-Ampt model is a relatively inferior approximation but in reality such information is often unavailable. We also concede that the choice of reference model still remains an open challenge for multi-model procedures and is beyond the scope of this paper. We hope to address such problem quantitively in future works. 

%

%%%%%%%%%%%%%%%%%%%%%%%%%%%%%%%%%%%%%%%%%%%%%%%%%%
\subsection{Effects of time-dependent model error}
\label{subsec:model-error-time}

\begin{figure}[htbp]
\begin{center}
\includegraphics[width=0.49\textwidth]{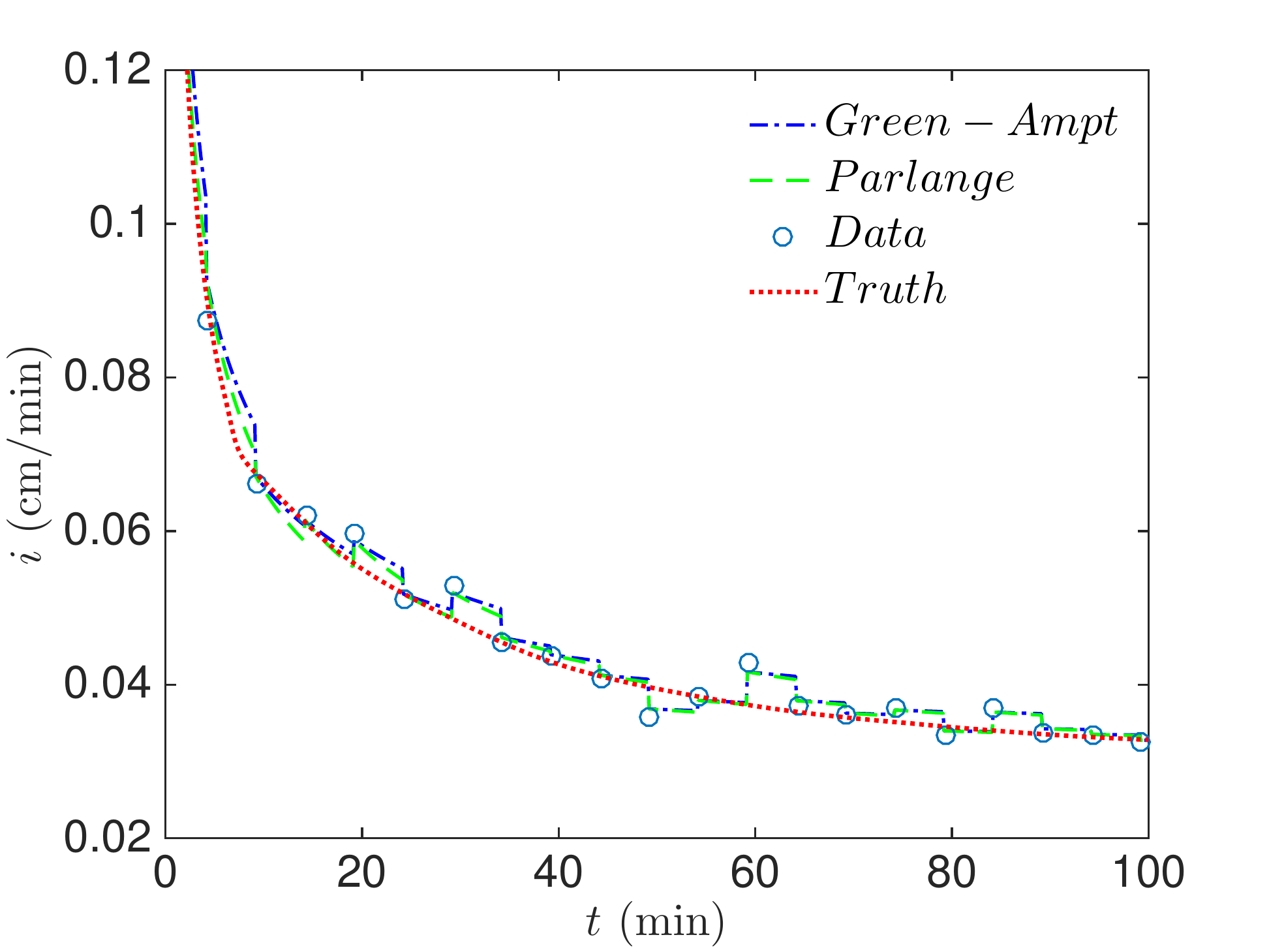} 
\includegraphics[width=0.49\textwidth]{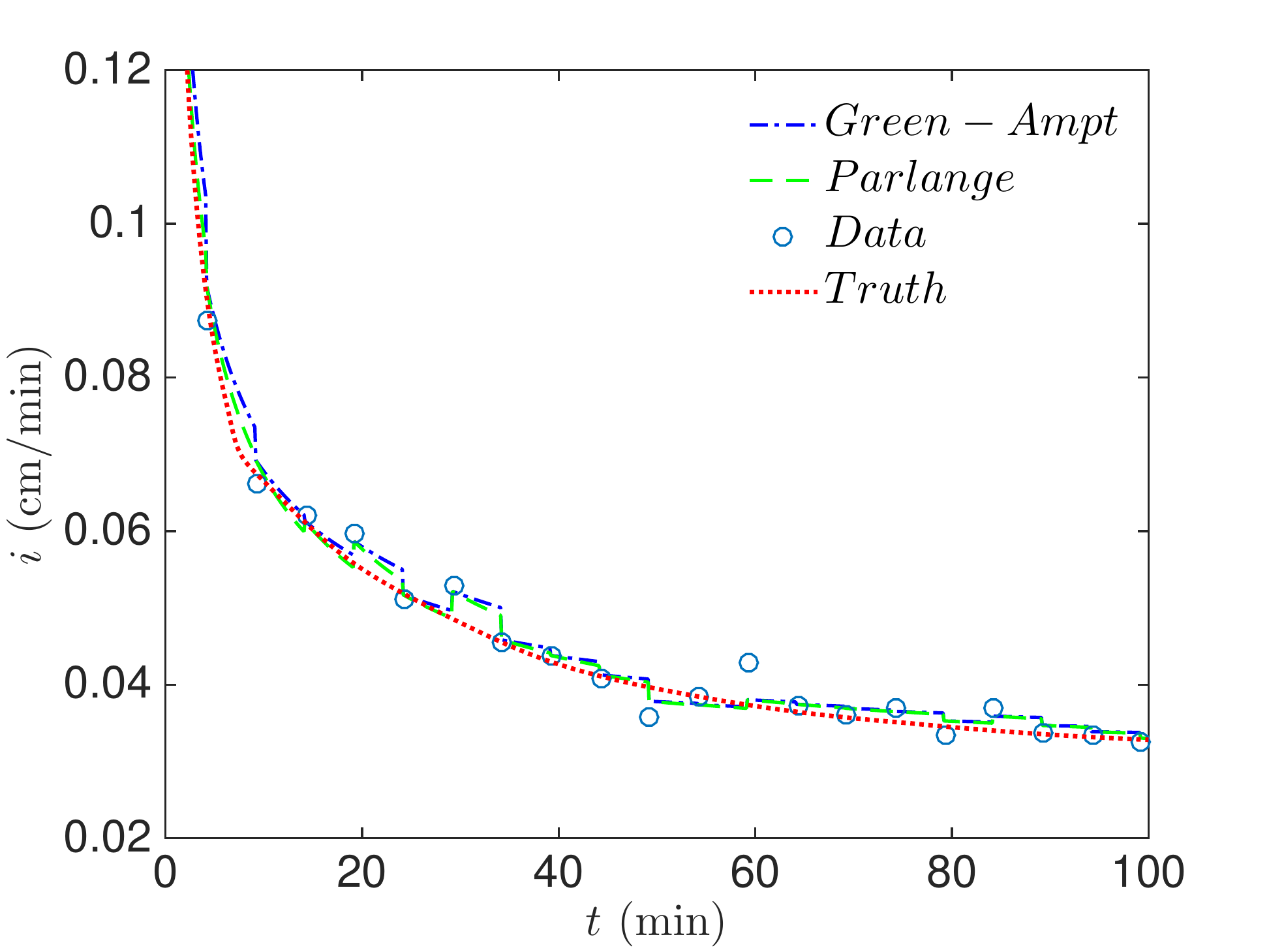} 
\includegraphics[width=0.49\textwidth]{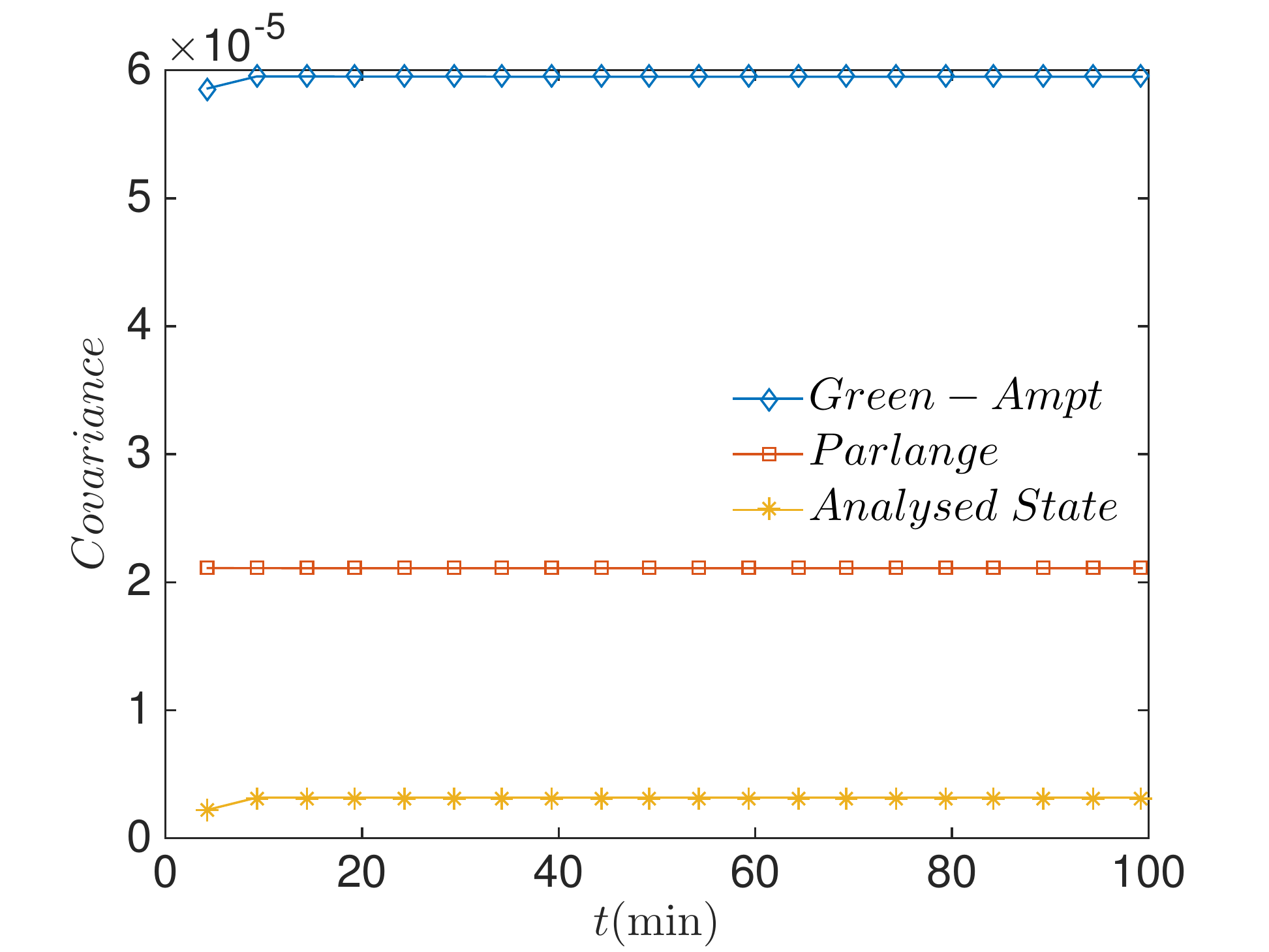} 
\includegraphics[width=0.49\textwidth]{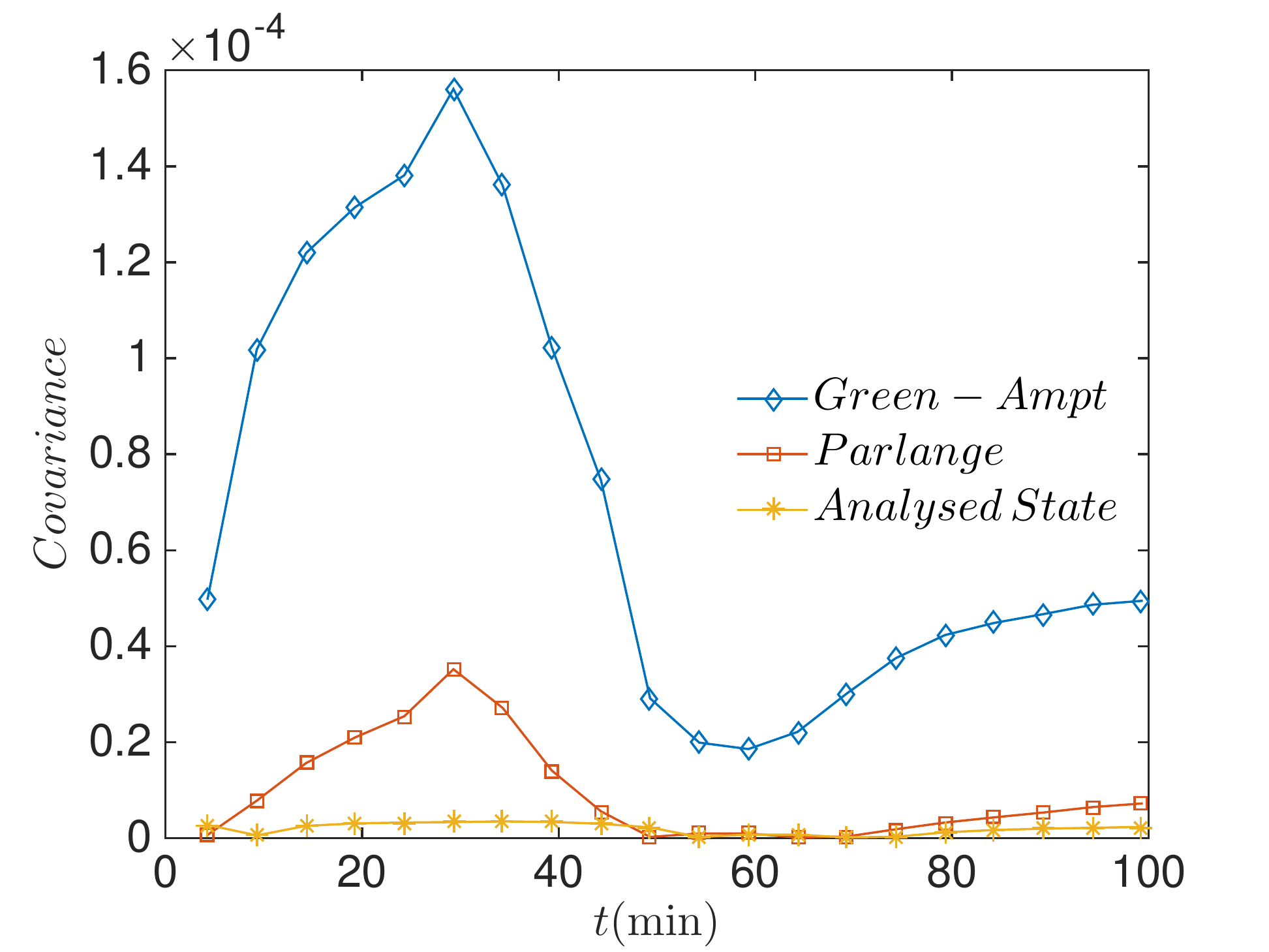} 
\includegraphics[width=0.49\textwidth]{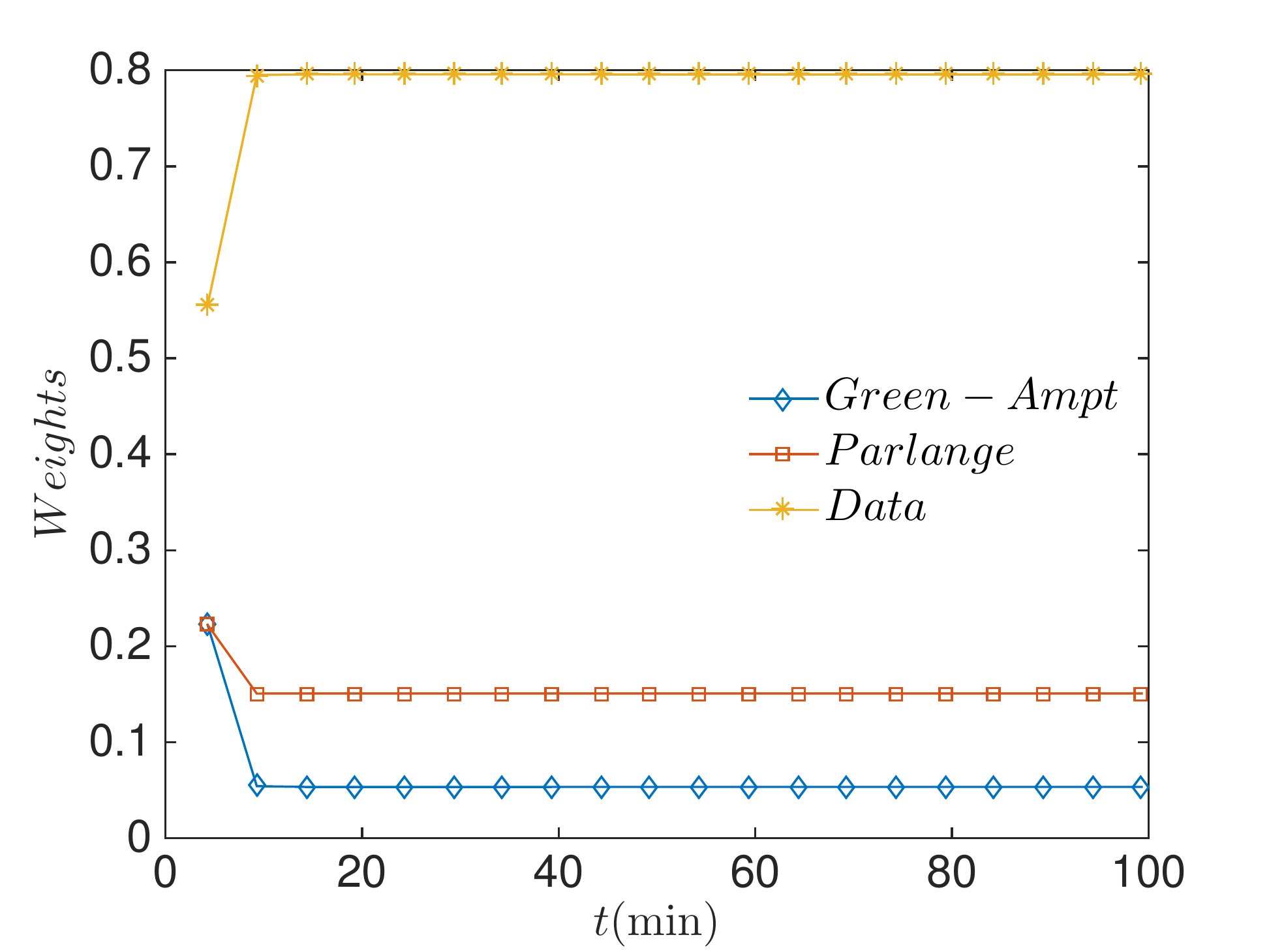} 
\includegraphics[width=0.49\textwidth]{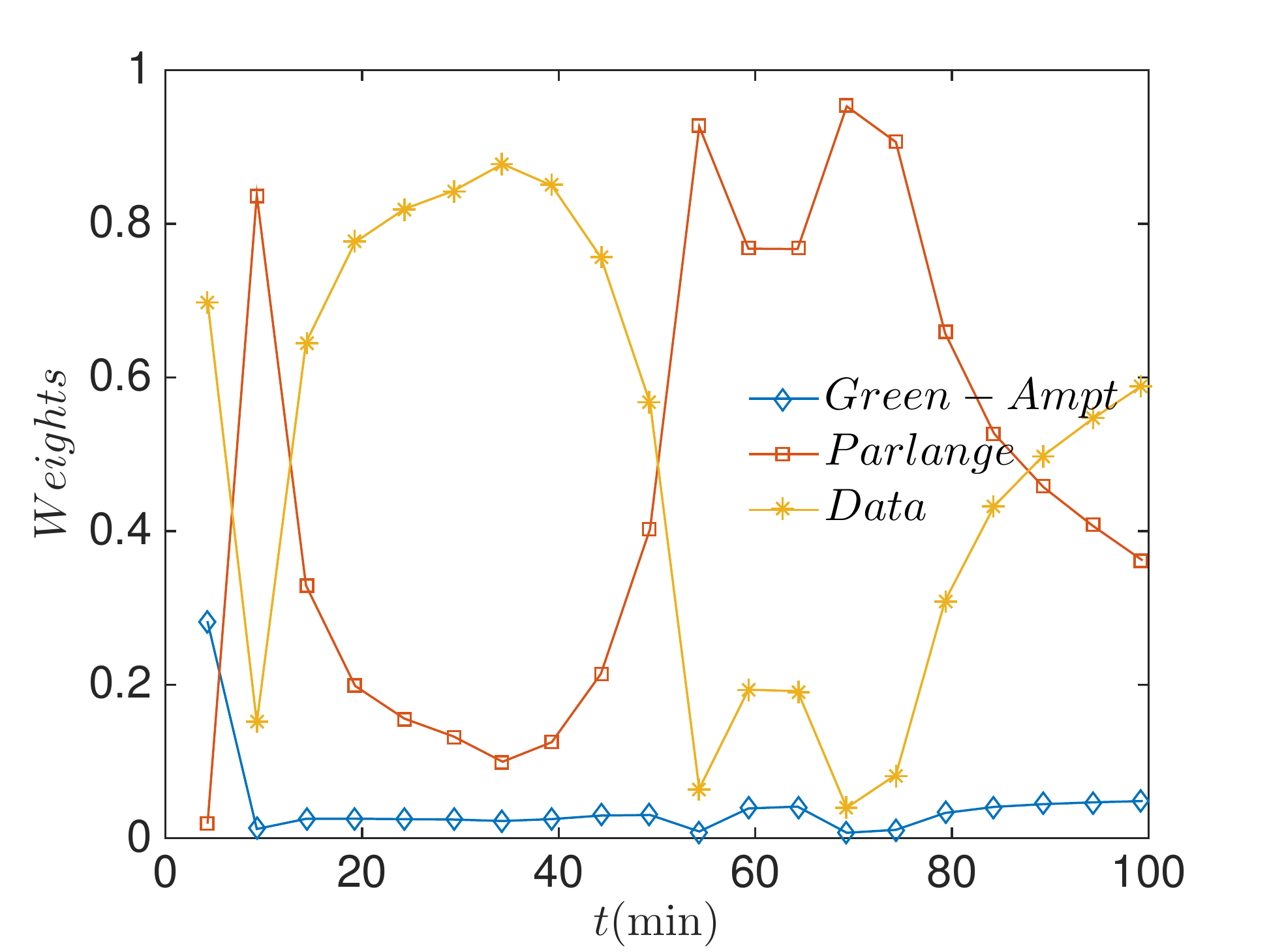} 
\caption{Temporal evolution of infiltration rate $i(t)$ assimilated from Green-Ampt model (dash-dot line), Parlange model (dashed line) and data (circle), using extended Kalman filter. The true value computed from VS2DT is denoted as dotted line. Model Covariances and weights are plotted as a function of time in the second and third rows, respectively. In the left plots, all models' error are taken as white noises; while they vary with time in the right plots.} 
\label{fig:model-error} 
\end{center}
\end{figure}

Data error is frequently due to the imperfectness of measuring instruments or human factors; hence it is frequently modelled as a time series of un-correlated random variables, e.g. a white noise process. However, factors that contribute to the discrepancy between system true dynamics and its mathematical models often have temporal structure and errors between different times are not uncorrelated. To investigate the impact of such time-dependence on assimilation results, we employ the EKF to conduct two numerical experiments in Fig. \ref{fig:model-error} using the same data set $\d$, but the model errors in the left column are homogenous random variables (whose variance is the mean error over time); while the right column's are heterogeneous (actual errors at each time step from Fig. \ref{fig:three-models}).  

From a glance at the temporal evolution of infiltration rates in Fig. \ref{fig:model-error}, EKF using time-dependent model errors produces more accurate assimilation results. Each model's predictive covariance $\U$ and weight converge to constants if the models errors are set as white noises. This is expected from earlier studies on the Kalman filter with a single model \cite{Grewal-1993-Kalman}. However, in case of time-heterogeneous $\e_m(t)$, each model's predictive covariance and weight flucuatate wildly with time and their values are notably different from their counterparts in the white noise setting. For example, in the left plots, the predictive covariance from the ``analyzed state" is considerably lower than the one from either Green-Ampt or Parlange model; and the highest weight is assigned to data, an indication that both model forecasts are ``worse" than the measurements. But in the right plots, the covariance from ``analyzed state" occasionally match with that from Parlange model, when the latter carries the heaviest weight ($t\in [50, 80]$). Thus the assimilation procedure could be improved by incorporating time-dependent model errors.

%

%%%%%%%%%%%%%%%%%%%%%%%%%%%%%%%%%%%%%%%%%%%%%%%%%%
\subsection{Comparison with Bayesian model averaging}
\label{subsec:BMA}

\begin{figure}[htbp]
\begin{center}
\includegraphics[width=0.49\textwidth]{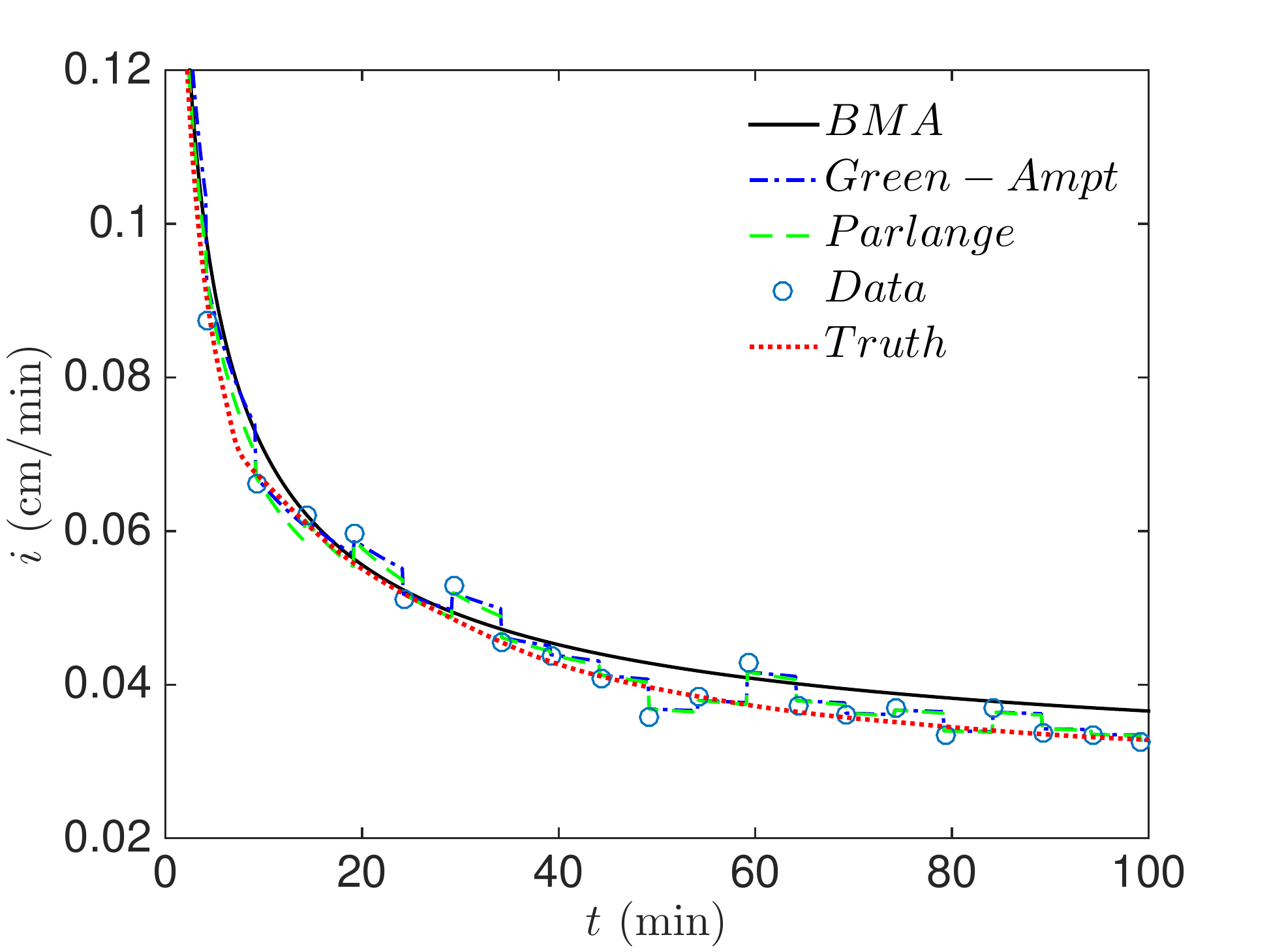} 
\includegraphics[width=0.49\textwidth]{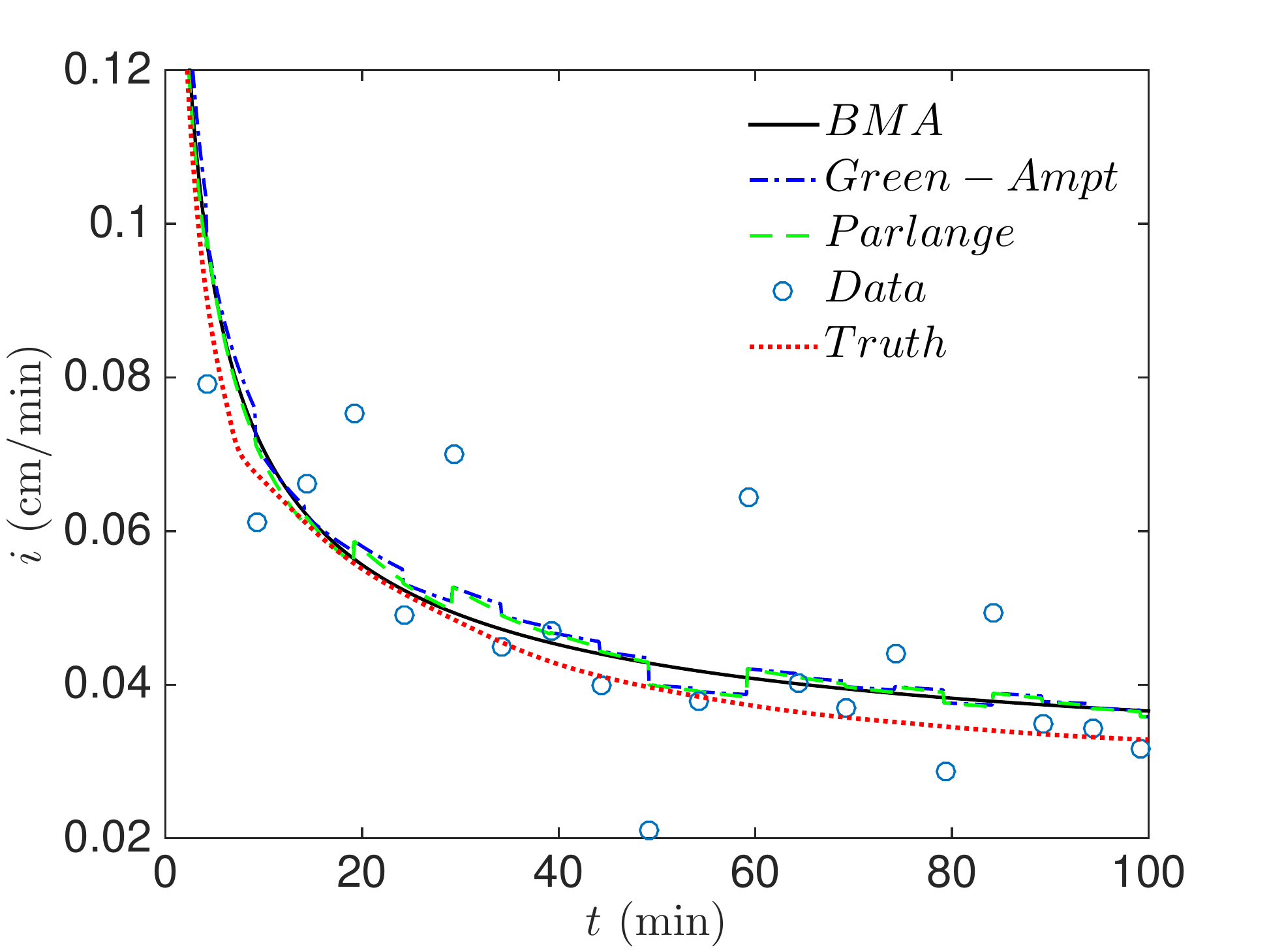} 
\caption{Temporal evolution of infiltration rate $i(t)$ assimilated from Green-Ampt model (Dash-dot line), Parlange model (dashed line) and data (circle), using extended Kalman filter; solutions from Bayesian model averaging is also presented (solid line). Noise of the measurement is set as $(1)$ $\e_d \sim \mathcal N (0,0.002^2)$ and $(2)$ $\e_d \sim \mathcal N (0,0.01^2)$. The true value computed from VS2DT is denoted as dotted line.} 
\label{fig:BMA-DA} 
\end{center}
\end{figure}

Bayesian model averaging (BMA) is a popular statistical procedure that infers predictions by weighing multiple model forecasts \cite{hoeting:1999:bma,raftery:2005:ubm,duan2007multi}. To evaluate its performance against the multi-model DA approach (and hence investigate the relative importance of data), we compute the temporal evolution of the infiltration rate from the EKF algorithm and BMA. Figure \ref{fig:BMA-DA} illustrates their prediction results using two data sets: one with smaller error (left plot), $\e_d \sim \mathcal N (0,0.002^2)$;  while the other's is noiser (right plot), $\e_d \sim \mathcal N (0,0.01^2)$. The BMA results do not change with respect to any data variation; but the DA approach is sensitive to the data input and a more accurate measurement would greatly improve its overall prediction.

\begin{figure}[htbp]
\begin{center}
\includegraphics[width=4in]{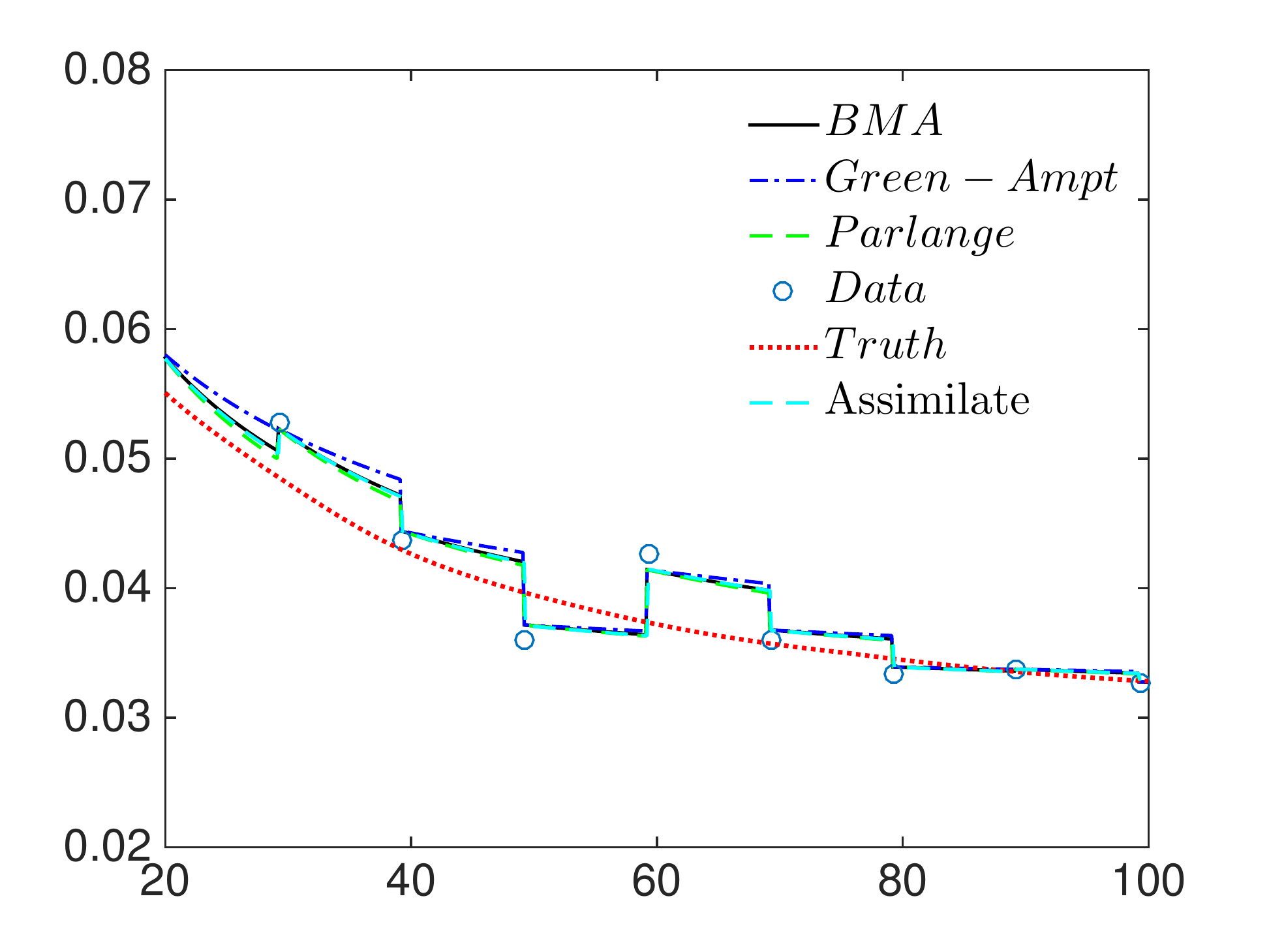} 
\caption{A segment of the infiltration rate $i(t)$ assimilated from extended Kalman filter using Green-Ampt model, Parlange model and data (circle). Its prediction from Bayesian model averaging (solid line) and Assimilated models (Green-Ampt $\&$ Parlange) (dashed line) is also presented. The true value computed from VS2DT is denoted as dotted line.} 
\label{fig:bma-data-short} 
\end{center}
\end{figure}

From the study above, one may conclude that BMA is structured to produce the most reliable prediction in the absence of data. Alternatively, \cite{Narayan-2012-Sequential} also states that the multi-model assimilation process is applicable when there are only model forecasts. We examine their merits by looking closer to the infiltration rate at time intervals between each observation. From Fig. \ref{fig:bma-data-short}, one can see that results from BMA and the assimilated Green-Ampt and Parlange models are identical, which confirms that our DA framework satisfies the same criteria as BMA in the absence of data.

%However, we note that discrepancy occurs in the period between the first two data points, which . The former finding confirms that , 

%while the former usually lie between the latter two forecasts. In other words, it is no better than the ``best" model forecasts. However, in practice, true system value is usually not a priori that complicates the models' assessments. BMA might not yield a better prediction but with the best of our knowledge.  

%

\section{Conclusion}\label{sec:summary}

In this paper we first present a novel framework for sequential data assimilation with multiple sources of models and data. Based on an earlier method resembling generalised Kalman filter, we employ extended Kalman Filter, ensemble Kalman filter and particle filter to investigate the impacts of model-form uncertainty. The particle filter extension we propose allows nonlinear models with non-Gaussian noise.  Our examples lead to the following empirical conclusions: 

\begin{itemize}

\item {Multi-model assimilation (extended Kalman filter, ensemble Kalman filter or particle filter) has the potential to offer a more accurate prediction than a single model forecast or data. }

%\item{Our particle filter data assimilation algorithm depends on choice of a ``reference model". The particular choice of which model is selected as the reference alters assimilation results. Thus, our procedure does not offer a constructive solution for inferring a hierarchy of models based on their errors.}

\item{Both the assimilation order and prior understanding of model accuracy matter in our multi-model particle data assimilation algorithm.}

\item{Time-dependent model errors can improve data assimilation results. How to model these time-dependent errors is still an open issue.}

\item {Our assimilation process yields identical results as Bayesian Model Averaging in the absence of observations. However, our algorithm superscedes BMA in applicability since it can also handle several sources of experimental or measurement data.}

%\item {Particle filter offers good improved predictions.}

\end{itemize}

%\input Appendix

%\section*{Acknowledgment}
%This work is in part supported by DARPPA, AFOSR, DOE/NNSA, NSF, CNSF.

\bibliographystyle{siamplain}
\bibliography{Bib/DA,Bib/random,Bib/collocation,Bib/bma-kalman,Bib/approximation,Bib/subflow}

\end{document}